\documentclass[11pt,a4paper]{scrartcl}

\usepackage[centertags]{amsmath} 
\usepackage{color} 
\usepackage{amsfonts}
\usepackage{amssymb}
\usepackage{amsthm}
\usepackage{epsfig} 
\usepackage{rotating}
\usepackage{psfrag}
\usepackage{caption}
\usepackage{booktabs}
\usepackage{enumitem}
\usepackage{pdflscape}

\newcommand{\R}{\mathbb{R}} 
\newcommand{\N}{\mathbb{N}} 
\newcommand{\C}{\mathbb{C}}

\newcommand{\X}{\mathcal{X}}

\newcommand{\D}{\mathcal{D}}
\newcommand{\Cc}{\mathcal{C}}
\newcommand{\Pp}{\mathcal{P}} 
 
\newcommand{\Rr}{\mathcal{R}} 
\newcommand{\Ss}{\mathcal{S}}
\newcommand{\W}{\mathcal{W}} 
\newcommand{\T}{\mathcal{T}} 
\newcommand{\V}{\mathcal{V}} 
\newcommand{\U}{\mathcal{U}} 
\newcommand{\Q}{\mathcal{Q}}

\newcommand{\lb}[1]{\underline{#1}} 
\newcommand{\ub}[1]{\overline{#1}} 
\newcommand{\blind}[1]{\textcolor{white}{#1}}

\newtheorem{thm}{Theorem}
\newtheorem{cor}[thm]{Corollary}
\newtheorem{lem}[thm]{Lemma}
\newtheorem{assum}{Assumption}

\newtheorem{defn}{Definition}
\newtheorem{exmp}{Example}
\newtheorem{rem}{Remark}
\newtheorem{alg}{Algorithm}

\begin{document}

\begin{center}
% Title
{\LARGE
\bf
\textsf{How scaling of the disturbance set affects\\[2mm]  robust positively invariant sets for linear systems}
}

\renewcommand{\thefootnote}{$\dagger$} 

\vspace{5mm}
{
Moritz Schulze Darup\footnotemark[1], Rainer Manuel Schaich\footnotemark[1], and Mark Cannon\footnotemark[1]
}
\vspace{2mm}

  \footnotetext[1]{M. Schulze Darup, R. M. Schaich, and M. Cannon are with 
the Control Group, Department of Engineering Science,
        University of Oxford, Parks Road, Oxford OX1 3PJ, UK.
        E-mail: {\tt moritz.schulzedarup@rub.de}.}
\end{center}

\paragraph{Abstract.}
This paper presents new results on robust positively
invariant (RPI) sets for linear discrete-time systems
with additive disturbances. In particular, we study how
RPI sets change with scaling of the disturbance set. 
More precisely, we show that many properties of RPI sets crucially depend on a unique scaling factor which determines the transition from nonempty to empty RPI sets. We characterize this critical scaling factor, present an efficient algorithm for its computation, and analyze it for a number of examples from the literature.


\section{Introduction}

Robust positively invariant (RPI) sets are important for performance analysis and synthesis of controllers for uncertain systems (see, e.g.,  \cite[Sects.~6.4 and 6.5]{Blanchini1999} or \cite[Sect.~VII]{Glover1971}).
In particular, RPI sets can be used to design robust model predictive control (MPC) schemes with guaranteed stability (see, e.g., \cite{Kouvaritakis2015,Lee1999,Mayne2006,Mayne2005,Rakovic2012}). In this paper, we address RPI sets for linear disturbed systems 
\begin{equation}
\label{eq:system}
x(k+1) = A \, x(k) + E \, d(k)
\end{equation}
with state and disturbance constraints of the form
\begin{equation}
\label{eq:constraints}
x(k) \in \X  \quad \text{and} \quad d(k) \in \D^\alpha \quad \text{for every} \quad k \in \N,
\end{equation}
where the set $\D^\alpha := \alpha \, \D^\ast$ denotes a scaled version  of a nominal disturbance set $\D^\ast$ for some scalar $\alpha>0$.
Roughly speaking, an RPI set $\Pp$ for system~\eqref{eq:system} with constraints~\eqref{eq:constraints} is such that the trajectory of the disturbed system~\eqref{eq:system} remains in $\Pp$ at all times $k\in \N$ for every initial condition $x_0 \in \Pp$ and for all disturbances $d(k) \in \D^\alpha$ (see Def.~\ref{def:RPI} further below for a formal definition of RPI sets).
Now, in this study, we are particularly interested how RPI sets for system~\eqref{eq:system} with constraints~\eqref{eq:constraints} change with the scaling $\alpha$. As we obviously have $\D^{\alpha_1} \subset \D^{\alpha_2} $ for $\alpha_1 < \alpha_2$, it is easy to see that the size of the maximal robust positively invariant (MRPI) set for \eqref{eq:system} with \eqref{eq:constraints} (see Def.~\ref{def:MRPI}) is non-increasing with~$\alpha$.
In addition, it is intuitively clear that the MRPI may be equal to the empty set for large scalings $\alpha$ (e.g., if $\alpha$ is such that $E \D^\alpha \nsubseteq \X$). In this paper, we show that there always exists a unique scaling factor $\alpha^\ast$, which defines the transition from a nonempty to an empty MRPI set. More precisely, we prove that (if $\X$ is bounded) there always exists an $\alpha^\ast>0$ such that the MRPI set is nonempty for every $\alpha \in (0,\alpha^\ast]$ but empty for every $\alpha > \alpha^\ast$. 

In principle, the existence of such a critical scaling factor (CSF)  $\alpha^\ast$ is neither surprising nor new. Similar observations were made, for instance, for the numerical examples discussed in 
\cite[Exmp.~4.1]{Kolmanovsky1995}, \cite[p. 114]{Kouvaritakis2015}, and \cite[Sect.~4]{Lee1999}. 
Moreover, we previously analyzed linear systems with scaled disturbance sets in the context of parametric RPI sets in \cite{Schaich2015} and \cite{SchulzeDarup2016_CDC}. 
In this paper, however, we analyze the CSF and its influence on properties of RPI sets in a more general framework than in 
\cite{Kolmanovsky1995,Kouvaritakis2015,Lee1999,Schaich2015,SchulzeDarup2016_CDC}.
First, and most importantly, we introduce a number of novel properties of RPI sets as summarized in Thm.~\ref{thm:criticalScaling}. In particular, we show that the CSF has crucial impact on the continuity, finite determinedness, and shape of RPI sets (see statements (iii)--(vii) in Thm.~\ref{thm:criticalScaling}).
Second, in contrast to \cite{Kolmanovsky1995,Kouvaritakis2015,Lee1999}, the analysis in this paper is not limited to specific examples.
Third, in contrast to \cite{Schaich2015,SchulzeDarup2016_CDC},
the results obtained  do not require $\X$ and $\D^\ast$ to be polytopic or $E\,\D^\ast$ to be full-dimensional in $\R^n$.
 Finally, we provide an efficient algorithm to numerically over- and underestimate the CSF with arbitrary precision (see Thm.~\ref{thm:algorithm}), which enables the proposed analysis scheme to be used. It is important to note that the evaluation of the CSF is useful even if the disturbance set is fixed for a given system (e.g., if $\alpha=1$). In fact, knowledge of the CSF can be of interest for determining allowable state and input constraints and acceptable strengths of disturbances when designing actuators, sensors, and controllers. In other words, the results in this paper are not only relevant for systems with disturbances that are contained in sets with variable scalings as in~\eqref{eq:constraints} but also for the classical case where the disturbances satisfy $d(k)\in \D^\ast$ (for every $k\in \N$).

The paper is organized as follows. In Sect.~\ref{sec:NotationDefinitions}, we collect definitions, assumptions, and known characteristics of RPI sets. The main result of the paper, i.e., how scaling of the disturbance set affects the properties of RPI sets %, in particular whether non-empty RPI sets exist, 
is discussed in Sect.~\ref{sec:effectScaling}. 
In this context, we also show that the effect of scaled disturbances  on RPI sets is different from the effect of scaled input constraints on controlled invariant (CI) sets (which are structurally related to RPI sets).
Then,  in Sect.~\ref{sec:algorithm}, an algorithm for the efficient computation of the CSF is presented. The algorithm is applied to a number of numerical examples in Sect.~\ref{sec:examples}.  Finally, conclusions are stated in Sect.~\ref{sec:conclusions}.

\section{Definitions, Assumptions, and Preliminaries}
\vspace{-2pt}

\label{sec:NotationDefinitions}

\subsection{Notation and set operations}

We denote positive real and natural numbers by $\R_+$ and $\N_+$, respectively. The notation $\N_{[i,k]}$ refers to $\N_{[i,k]}:=\{ j \in \N \,|\, i \leq j \leq k \}$.
Now, let $p \in \N_+$ and consider two bounded and convex sets $\U , \V \subset \R^p$ containing the origin. We frequently use the following manipulations of sets.
 The scaling of a set by some factor $\beta \in \R_+$  is defined as $\beta \,\U := \{ \beta u \in \R^p \,|\, u \in \U\} $.
Moreover, for $q \in \N_+$ and some matrices $C \in \R^{p\times p}$ and $D \in \R^{q \times p}$, we define $C^{-1} \, \U:= \{ \xi \in \R^p \,|\, C\,\xi \in \U \}$ and  $D\,\U := \{ D u \in \R^q \,|\, u \in \U \} $. Note that $C^{-1} \, \U$ is well defined even if $C$ is not invertible. Finally, the operations
\begin{align*}
\U \oplus \V &:= \{ u+v \in \R^p \,|\,  u \in \U, v \in \V\} \qquad \text{and} \\
\U \ominus \V &:= \{ \xi \in \R^p \,|\, \forall v \in \V: \xi+v \in \U\},
\end{align*}
describe the Minkowski addition and the Pontryagin difference, respectively. 
It is easy to see that both operations as well as the intersection of two sets are distributive in the sense that 
\begin{equation}
\label{eq:distributiveLaws}
\beta \,( \U \oplus \V) = (\beta \U) \oplus (\beta \V), \quad  \beta \,( \U \ominus \V) = (\beta \U) \ominus (\beta \V),  \quad \text{and} \quad \beta \,( \U \cap \V) = (\beta \U) \cap (\beta \V)
\end{equation}
for every $\beta \in \R_+$. 
%Eventually, we use the following shorthand notations.
A C-set is a convex and compact set containing the origin as an interior point. The boundary, the interior, and the closure of a set $\U$ are denoted by $\partial \U$, $\textrm{int}(\U)$, and $\textrm{cl}(\U)$, respectively. 
The matrix $I_p$ refers to the identity matrix in $\R^{p \times p}$. %The vector $0_p$ denotes the origin in $\R^p$. 
The rank of a matrix $D$ is denoted by $\textrm{rk}(D)$.
A $p$-dimensional vector with with all entries equal to $1$ is written as $\boldsymbol{1}_p$.

\subsection{Formal definition of robust positively invariant sets}

The following statements provide a precise definition of robust positive invariance. In every definition, we assume $\X\subseteq \R^n$, $\D^\ast \subset \R^m$, and $\alpha \in \R_+$.

\begin{defn}
\label{def:RPI}
A set $\Pp \subseteq \R^n$ is called robust positively invariant (RPI) for system~\eqref{eq:system} with constraints~\eqref{eq:constraints}, if (i) $\Pp \subseteq \X$ and (ii) $A\, x+E\,d\in \Pp$ for every $x \in \Pp$ and every $d \in \D^\alpha$.
\end{defn}

Note that both the union and the intersection of two RPI sets result in another RPI set. Further note that $\Pp=\emptyset$ is a RPI set according to Def.~\ref{def:RPI}. Based on these observations, the following definitions of the maximal and minimal RPI set are reasonable. 

\begin{defn}
\label{def:MRPI}
The union of all RPI sets for~\eqref{eq:system} with~\eqref{eq:constraints} is called the maximal robust positively invariant (MRPI) set for~\eqref{eq:system} with~\eqref{eq:constraints} and denoted by $\Pp_{\max}^\alpha$.
\end{defn}

\begin{defn}
\label{def:mRPI}
The intersection of all nonempty RPI sets and the MRPI set for~\eqref{eq:system} with~\eqref{eq:constraints} is called the minimal robust positively invariant (mRPI) set for~\eqref{eq:system} with~\eqref{eq:constraints} and denoted by $\Pp_{\min}^\alpha$.
\end{defn}

To fully understand Def.~\ref{def:mRPI}, note that, if nonempty RPI sets exist, the mRPI denotes the ``smallest'' nonempty RPI set. Clearly, in this case, the intersection with the MRPI set is irrelevant since the MRPI is one of the nonempty RPI sets.
If, however, no nonempty RPI set exists, the intersection with the MRPI set (which will, in this case, be the empty set) avoids an empty intersection (which would be unequal to the empty set).

\subsection{Assumptions on the system matrices and constraints}

Most statements in the paper require the two following assumptions.
Note that the structure in~\eqref{eq:AE} can be derived for every (stabilizable) pair $(A,E)$ %with $\mathrm{rk}(E) \geq 1$ 
using a suitable linear transformation.

\begin{assum}
\label{assum:matricesAndSets}
The system matrices $A \in \R^{n \times n}$ and $E \in \R^{n \times m}$ are structured as 
\begin{equation}
\label{eq:AE}
A = \begin{pmatrix}
A_{11} & A_{12} \\ 0 & A_{22}
\end{pmatrix} \qquad \text{and} \qquad E = \begin{pmatrix}
E_{1}  \\ 0 
\end{pmatrix}
\end{equation}
with $A_{11} \in \R^{r \times r}$, $A_{12} \in \R^{r \times (n-r)}$, $A_{22} \in \R^{(n-r) \times (n-r)}$, and $E_1 \in \R^{r \times m}$, where $r \in \N_{[1,n]}$. The pair $(A_{11},E_1)$ is controllable. The sets $\X \subset \R^n$ and $\D^\ast \subset \R^m$ are C-sets and $\alpha \in \R_+$.
\end{assum}

\begin{assum}
\label{assum:stability}
The eigenvalues $\lambda \in \C$ of the matrices $A_{11}$ and $A_{22}$ in Assum.~\ref{assum:matricesAndSets} are strictly stable (i.e., $| \lambda|<1$).
\end{assum}

\subsection{Known properties of the maximal and minimal robust positively invariant sets}

Definitions~\ref{def:MRPI} and \ref{def:mRPI} immediately lead to the following relation between the MRPI and the mRPI set.

\begin{lem}
\label{lem:PminPmaxEmpty}
Let Assums.~\ref{assum:matricesAndSets} and~\ref{assum:stability} be satisfied. Then, (i) $\Pp_{\min}^\alpha$ is empty if and only if $\Pp_{\max}^\alpha$ is empty and (ii) $\Pp_{\min}^\alpha \subseteq \Pp_{\max}^\alpha \subseteq \X$.
\end{lem}

In order to characterize $\Pp_{\max}^\alpha$ and $\Pp_{\min}^\alpha$ more precisely, we will analyze two special sequences of sets (inspired by the sequences in \cite[Eqs. (1.10) and (4.3)]{Kolmanovsky1998}).
The first sequence 
\begin{equation}
\label{eq:sequenceSk}
\Ss_{k+1}^\alpha := A^{-1} (\Ss_k^\alpha \ominus E \, \D^\alpha) \cap \X \qquad \text{with} \qquad  \Ss_{0}^\alpha:=\X,
\end{equation}
provides sets $\Ss_k^\alpha$ that contain initial conditions $x(0)$ for which the trajectory of the disturbed system~\eqref{eq:system} remains in $\X$ for at least $k$ time steps and for all disturbances $d(j) \in \D^\alpha$.
The second sequence
\begin{equation}
\label{eq:sequenceRk}
\Rr_{k+1}^\alpha := A \Rr_k^\alpha \oplus E \D^\alpha \qquad \text{with} \qquad  \Rr_{0}^\alpha:=\{0\}.
\end{equation}
defines sets $\Rr_k^\alpha$ that contain states $x^\ast$ for which there exist disturbance sequences $d(0),\dots,d(k-1)$  $\in \D^\alpha$ such that
$x^\ast = \sum_{j=0}^{k-1} A^{k-1-j} E \,d(j)$.
In other words, $\Rr_k^\alpha$ contains states that are reachable  from the origin in $k$ steps. As summarized in the following lemmas and remark, the  sequences~\eqref{eq:sequenceSk} and~\eqref{eq:sequenceRk} have well-known properties (see, e.g., \cite{Kolmanovsky1998,Hirata2003})

\begin{lem}
\label{lem:SkAlphaSLimit}
Let Assums.~\ref{assum:matricesAndSets} and~\ref{assum:stability} be satisfied. Then, the following statements hold.
\begin{enumerate}
\item[(i)] For every $k \in \N$, the set $\Ss_k^\alpha$ is either a C-set in $\R^n$ or empty.
\item[(ii)] The sequence~\eqref{eq:sequenceSk} is non-increasing, i.e.,  $\Ss_{k+1}^\alpha \subseteq \Ss_k^\alpha$ 
 for every $k \in \N$.
\item[(iii)] The limit of the sequence~\eqref{eq:sequenceSk} is
$\Ss_\infty^\alpha:= \bigcap_{k=0}^\infty \Ss_k^\alpha$.
\item[(iv)] The limit set $\Ss_\infty^\alpha$ is either compact and convex with $0 \in \Ss_\infty^\alpha$ or empty.
\end{enumerate}
\end{lem}

\begin{lem}
\label{lem:RkAlphaRLimit}
Let Assums.~\ref{assum:matricesAndSets} and~\ref{assum:stability} be satisfied. Then, the following statements hold.
\begin{enumerate}
\item[(i)] For every $k\in \N$, the set $\Rr_k^\alpha$ is compact and convex with $0 \in \Rr_k^\alpha$. 
\item[(ii)] The sequence~\eqref{eq:sequenceRk} is non-decreasing, i.e.,  $\Rr_{k+1}^\alpha \supseteq \Rr_k^\alpha$ 
 for every $k \in \N$.
\item[(iii)] The limit set of the sequence~\eqref{eq:sequenceRk} is
$\Rr_\infty^\alpha:=   \bigcup_{k=0}^\infty \Rr_k^\alpha$.
\item[(iv)] The limit set $\Rr_\infty^\alpha$ is bounded and convex with $0 \in \Rr_\infty^\alpha$.
\end{enumerate}
\end{lem}

\begin{rem}
\label{rem:subsetDimr}
The structure of~\eqref{eq:AE} implies certain structural properties of the sets $\Rr_{k}^\alpha$.
In fact, defining the projection matrices $P_1:= \begin{pmatrix}
I_r & 0
\end{pmatrix} \in \R^{ r \times n}$ and $P_2 := \begin{pmatrix}
0 & I_{n-r} 
\end{pmatrix} \in \R^{(n-r) \times n}$ and assuming $r<n$, it is easy to see that $P_2 \, \Rr_{k}^\alpha = \{0\}$ for every $k \in \N$. Moreover, $P_1 \, \Rr_{k}^\alpha$ is a C-set in $\R^r$ for every $k \geq r$ (cf.~\cite[Rem.~3]{Hirata2003}). Thus, if one is interested in the explicit computation of the sets $\Rr_k^\alpha$, it would make sense to only consider the nontrivial $r$-dimensional subspace. However, for the derivation of the theoretical results in this paper, it is more convenient to address $\Ss_k^\alpha$ and $\Rr_k^\alpha$ in the same $n$-dimensional space.
\end{rem}

It is well-known that the limits $\Ss_\infty^\alpha$ and $\Rr_\infty^\alpha$ are closely related to $\Pp_{\max}^\alpha$ and $\Pp_{\min}^\alpha$, respectively. In fact, according to Lem.~\ref{lem:PmaxPmin}, $\Ss_\infty^\alpha$ is the MRPI set and, if $\Rr_{\infty}^\alpha\subseteq \X$, $\Rr_\infty^\alpha$ is the mRPI set.

\begin{lem}
\label{lem:PmaxPmin}
Let Assums.~\ref{assum:matricesAndSets} and~\ref{assum:stability} be satisfied. Then, $\Pp_{\max}^\alpha=\Ss_\infty^\alpha$ and 
\begin{equation}
\label{eq:Pmin}
\Pp_{\min}^\alpha=\left\{ \begin{array}{ll}
\Rr_{\infty}^\alpha & \text{if} \quad \Rr_{\infty}^\alpha\subseteq \X, \\
\emptyset & \text{otherwise}.
\end{array} \right.
\end{equation}
\end{lem}

Finally, while both sequences~\eqref{eq:sequenceSk} and~\eqref{eq:sequenceRk} are linked to RPI sets, they do not appear to be closely related to one another.
This conclusion is wrong, though, as the following lemma shows.
\begin{lem}
\label{lem:SkRk}
Let Assums.~\ref{assum:matricesAndSets} and~\ref{assum:stability} be satisfied and let $k \in \N$. Then, 
\begin{equation}
\label{eq:relationSkRk}
\Ss_{k}^\alpha = \bigcap_{j=0}^k (A^{j})^{-1} (\X \ominus \Rr_{j}^\alpha).
\end{equation}
\end{lem}
In principle, relation~\eqref{eq:relationSkRk} immediately follows from \cite[Eqs.~(5.1) and (5.2)]{Kolmanovsky1998}.
However, since our setup slightly differs from the one in \cite{Kolmanovsky1998}, we provide a formal proof in the appendix for completeness. 
Finally, from~\eqref{eq:relationSkRk}, it is easy to see that the relation 
\begin{equation}
\label{eq:relationSkRkPlus1}
\Ss_{k+1}^\alpha = (A^{k+1})^{-1} (\X \ominus \Rr_{k+1}^\alpha) \cap \Ss_k^\alpha
\end{equation}
holds for every $k \in \N$.

\section{The effect of scaled disturbances}
\label{sec:effectScaling}

In this section, we study the effect of variations in the scaling $\alpha$ on the properties of the MRPI and mRPI sets. In particular, we show that there always exists a scaling factor $\alpha^\ast>0$ such that  $\Pp_{\max}^\alpha$ and $\Pp_{\min}^\alpha$ are nonempty for $\alpha \in (0,\alpha^\ast]$ but empty for $\alpha > \alpha^\ast$.
While this observation is quite intuitive, we prove that (i) variations of $\Pp_{\max}^\alpha$ and $\Pp_{\min}^\alpha$ resulting from small changes in $\alpha$  and (ii) the finite determinideness of  $\Pp_{\max}^\alpha$ crucially depend on the CSF $\alpha^\ast$ (see Thm.~\ref{thm:criticalScaling} in Sect.~\ref{subsec:criticalStudyOfRPI}). 
In other words, many elementary properties of the MRPI and mRPI sets instantaneously change for scaling factors around $\alpha^\ast$.
We finally illustrate that this behavior is, in some sense, unique to RPI sets. In fact, while controlled invariant (CI) sets (for systems with controllable inputs) show some similarities to RPI sets, 
CI sets are not sensitive to changes in the scaling $\alpha$ (see Thm.~\ref{thm:MCIespDelta} in Sect.~\ref{subsec:similarStudyOfCI}).

\subsection{A critical study of robust positively invariant sets}
\label{subsec:criticalStudyOfRPI}

We begin with a more precise characterization of the CSF $\alpha^\ast$.
Obviously, there always exists an $\alpha^\ast \in \R_+$ such that $\Rr_\infty^\alpha \subseteq \X$ for every $\alpha \in (0,\alpha^\ast]$ but $\Rr_\infty^\alpha \nsubseteq \X$ for every $\alpha > \alpha^\ast$. According to Lems.~\ref{lem:PminPmaxEmpty} and \ref{lem:PmaxPmin}, such an $\alpha^\ast$ immediately implies that $\Pp_{\max}^\alpha$ and $\Pp_{\min}^\alpha$ are nonempty for $\alpha \in (0,\alpha^\ast]$ but empty for $\alpha > \alpha^\ast$.
It thus makes sense to formally define the CSF as
\begin{equation}
\label{eq:criticalAlpha}
\alpha^\ast :=  \sup \,\{\alpha \in \R_+ \,|\, \Rr^\alpha_\infty \subseteq \X \}.
\end{equation}
Note that the supremum is used in~\eqref{eq:criticalAlpha} since $\Rr^\alpha_\infty$ may or may not be closed.
%Further note that $\alpha^\ast$ is well-defined and finite 
We next show that $\alpha^\ast$ not only marks the transition from nonempty to empty sets $\Pp_{\max}^\alpha$ and $\Pp_{\min}^\alpha$.
In fact, as specified in statement (iii) of Thm.~\ref{thm:criticalScaling}, $\alpha^\ast$ also marks a ``discontinuity'' of the sets $\Pp_{\max}^\alpha$ and $\Pp_{\min}^\alpha$ as a function of $\alpha$.
Moreover, $\alpha^\ast$ plays an important role for the finite determinedness of $\Pp_{\max}^\alpha$ (see statement (iv)).
Finally, $\alpha^\ast$ is crucial for ``contact points'' between the boundaries of the sets $\Pp_{\min}^\alpha$, $\Pp_{\max}^\alpha$, and $\X$ (see statements (v) and (vi)).

\begin{thm}
\label{thm:criticalScaling}
Let Assums.~\ref{assum:matricesAndSets} and~\ref{assum:stability} be satisfied and let $\alpha^\ast$ be defined as in~\eqref{eq:criticalAlpha}. Then, the following statements hold.
\begin{enumerate}
\item[(i)] $\alpha^\ast$ is well-defined and finite.
\item[(ii)]  (a) $\Pp_{\min}^\alpha=\emptyset$ and (b) $\Pp_{\max}^\alpha=\emptyset$ if and only if $\alpha >\alpha^\ast$.
\item[(iii)] For every $\epsilon\in \R_+$ there exists a $\delta \in \R_+$ such that
\begin{subequations}
\begin{align}
\label{eq:PminPmaxEpsDeltaA}
\Pp_{\min}^{\alpha} &\subseteq  \Pp_{\min}^{\alpha+\delta} \subseteq (1+\epsilon) \,\Pp_{\min}^{\alpha}
\qquad \text{and}\\
\label{eq:PminPmaxEpsDeltaB}
 \Pp_{\max}^{\alpha+\delta} &\subseteq \Pp_{\max}^{\alpha} \subseteq (1+\epsilon) \,\Pp_{\max}^{\alpha+\delta} 
\end{align}
\end{subequations}
if and only if  $\alpha \neq \alpha^\ast$.
\item[(iv)] There exists a finite $k \in \N$ such that $\Pp_{\max}^\alpha = \Ss_k^\alpha$ if $\alpha \neq \alpha^\ast$.
\item[(v)] $\Pp_{\max}^\alpha$ is a C-set in $\R^n$ if $\alpha < \alpha^\ast$.
\item[(vi)]  (a) $\partial \Pp_{\min}^\alpha \cap \partial\X \neq \emptyset$ and (b) $\partial \Pp_{\min}^\alpha \cap \partial\Pp_{\max}^\alpha \neq \emptyset$  if and only if $\alpha = \alpha^\ast$.
\item[(vii)] $\partial \Pp_{\max}^\alpha \cap \partial \X \neq \emptyset$ if and only if $\alpha \leq \alpha^\ast$.
\end{enumerate} 
\end{thm}

The proof of Thm.~\ref{thm:criticalScaling} makes use of the two following lemmas.  Lemma~\ref{lem:RLimitAlpha} (which we prove in the appendix for completeness) formalizes an observation that was made in~\cite[p. 208]{Blanchini2008}. Lemma~\ref{lem:SkAffectedByDelta} provides the key to prove relation~\eqref{eq:PminPmaxEpsDeltaB}.

\begin{lem}
\label{lem:RLimitAlpha}
Let Assums.~\ref{assum:matricesAndSets} and~\ref{assum:stability} be satisfied. Then $\Rr_\infty^\alpha = \alpha \,\Rr_\infty^1$.
\end{lem}

\begin{lem}
\label{lem:SkAffectedByDelta}
Let Assums.~\ref{assum:matricesAndSets} and~\ref{assum:stability} be satisfied, let $\eta \in (0,1)$, and let $\alpha^\ast$ be defined as in~\eqref{eq:criticalAlpha}. Assume   $\alpha < \alpha^\ast$ and assume $\delta \in \R_+$ is such that $\delta \leq (1-\eta) (\alpha^\ast - \alpha)$. Then
\begin{equation}
\label{eq:etaSkAlphaSubsetSkAlphaPlusDelta}
\eta \, \Ss_k^\alpha  \subseteq \Ss_k^{\alpha+\delta} \subseteq \Ss_k^{\alpha}
\end{equation}
for every $k \in \N$.
\end{lem}
\begin{proof}[Proof of Lem.~\ref{lem:SkAffectedByDelta}]
We prove the claim by induction. Obviously, \eqref{eq:etaSkAlphaSubsetSkAlphaPlusDelta} holds for $k=0$ since $\Ss_0^\alpha=\Ss_0^{\alpha+\delta}=\X$ and since $\X$ being a C-set implies $\eta\, \X \subset \X $.
 It remains to show that \eqref{eq:etaSkAlphaSubsetSkAlphaPlusDelta}  implies
\begin{equation}
\label{eq:etaSkPlusOneAlphaSubsetSkAlphaPlusDelta}
\eta \, \Ss_{k+1}^\alpha  \subseteq \Ss_{k+1}^{\alpha+\delta} \subseteq \Ss_{k+1}^{\alpha}.
\end{equation}
The second relation in~\eqref{eq:etaSkPlusOneAlphaSubsetSkAlphaPlusDelta} can be proved using~\eqref{eq:sequenceSk}. In fact, we obviously have $\Ss_k^{\alpha+\delta} \ominus E \, \D^{\alpha+\delta} \subseteq \Ss_k^\alpha \ominus E \, \D^\alpha $
due to $\Ss_k^{\alpha+\delta} \subseteq \Ss_k^{\alpha}$ and $\D^{\alpha+\delta}\supset \D^\alpha$, and hence
$$
\Ss_{k+1}^{\alpha+\delta} = A^{-1} (\Ss_k^{\alpha+\delta} \ominus E \, \D^{\alpha+\delta}) \cap \X \subseteq A^{-1} (\Ss_k^\alpha \ominus E \, \D^\alpha) \cap \X = \Ss_{k+1}^{\alpha}.
$$
To prove the first relation in~\eqref{eq:etaSkPlusOneAlphaSubsetSkAlphaPlusDelta}, first note that~\eqref{eq:relationSkRkPlus1} implies
$$
\eta\,\Ss_{k+1}^\alpha = \eta\,(A^{k+1})^{-1} (\X \ominus \Rr_{k+1}^\alpha) \cap \eta \,\X \quad \text{and} \quad \Ss_{k+1}^{\alpha+\delta} = (A^{k+1})^{-1} (\X \ominus \Rr_{k+1}^{\alpha+\delta}) \cap \X.$$
Thus, the first relation in~\eqref{eq:etaSkPlusOneAlphaSubsetSkAlphaPlusDelta}
holds if 
\begin{equation}
\label{eq:conditionEtaXRkPlus1}
\eta \, (\X \ominus \Rr_{k+1}^\alpha) \subseteq \X \ominus \Rr_{k+1}^{\alpha+\delta}.
\end{equation}
According to Lem.~\ref{lem:RLimitAlpha} and~\cite[Thm. 2.1]{Kolmanovsky1998}, the l.h.s.~in~\eqref{eq:conditionEtaXRkPlus1} can be rewritten as
\begin{equation}
\label{eq:lhsRewritten}
\eta \, (\X \ominus \Rr_{k+1}^\alpha)=\eta \, (\X \ominus \alpha\,\Rr_{k+1}^1) = (\X \ominus \alpha \,\Rr_{k+1}^1) \ominus (1-\eta)\,(\X \ominus \alpha \,\Rr_{k+1}^1).
\end{equation}
Analogously, the r.h.s.~in~\eqref{eq:conditionEtaXRkPlus1} evaluates to
\begin{equation}
\label{eq:rhsRewritten}
\X \ominus \Rr_{k+1}^{\alpha+\delta} = \X \ominus (\alpha+\delta)\,\Rr_{k+1}^{1}= (\X \ominus \alpha \,\Rr_{k+1}^1) \ominus \delta \,\Rr_{k+1}^1.
\end{equation}
From comparison of Eqs.~\eqref{eq:lhsRewritten} and \eqref{eq:rhsRewritten}, we infer that condition~\eqref{eq:conditionEtaXRkPlus1} holds if 
\begin{equation}
\label{eq:condition1MinusEtaXMinusRk}
 (1-\eta)\,(\X \ominus \alpha \,\Rr_{k+1}^1) \supseteq \delta \,\Rr_{k+1}^1.
\end{equation}
Finally, an inner approximation of the set on the l.h.s.~in~\eqref{eq:condition1MinusEtaXMinusRk} can be derived according to \cite[Rem. 2.1]{Kolmanovsky1998}. In fact, due to $\Rr_{k+1}^1  \subseteq  \Rr_\infty^1$ since $\alpha^\ast\, \Rr_\infty^1 \subseteq \X$ by~\eqref{eq:criticalAlpha} and Lem.~\ref{lem:RLimitAlpha}, we obtain
$$
\X \ominus \alpha \,\Rr_{k+1}^1  \supseteq \left(\frac{\alpha^\ast}{\alpha} - 1\right)\alpha \,\Rr_{k+1}^1 = (\alpha^\ast - \alpha)\,\Rr_{k+1}^1.
$$
Thus, \eqref{eq:condition1MinusEtaXMinusRk} holds due to 
$\delta \leq (1-\eta)(\alpha^\ast - \alpha)$.
\end{proof}

\begin{proof}[Proof of Thm.~\ref{thm:criticalScaling}]
We separately prove the seven statements in Thm.~\ref{thm:criticalScaling}.

Statement (i). 
First note that $\Rr_\infty^\alpha = \alpha \,\Rr_\infty^1$ according to Lem.~\ref{lem:RLimitAlpha}. Thus, $\alpha^\ast$ from~\eqref{eq:criticalAlpha} refers to the largest scaling $\alpha \in \R_+$ of the set $\Rr_\infty^1$ such that $\alpha \,\Rr_\infty^1 \subseteq \X$. This scaling is well-defined and finite since $\X$ is a C-set, since $\Rr_\infty^1$ is bounded, and since  $\Rr_\infty^1$ contains more points than just the origin (see Rem.~\ref{rem:subsetDimr}).

Statement (ii). 
We only prove statement (ii).(a) since (ii).(b) then immediately follows from Lem.~\ref{lem:PminPmaxEmpty}. Clearly, the definition of $\alpha^\ast$ in~\eqref{eq:criticalAlpha}   implies   $\Rr_\infty^\alpha \nsubseteq \X$
for every $\alpha> \alpha^\ast$.
According to  Lem.~\ref{lem:PmaxPmin}, we thus have  $\Pp_{\min}^\alpha = \emptyset$    whenever $\alpha> \alpha^\ast$.  
It remains to show that $\Pp_{\min}^\alpha \neq \emptyset$ for every $\alpha \in (0,\alpha^\ast]$. This, however, easily follows from Eq.~\eqref{eq:criticalAlpha} in combination with Lems.~\ref{lem:PmaxPmin} and \ref{lem:RLimitAlpha}.

Statement (iii). First note that, for every $\alpha>\alpha^\ast$ and every $\delta \in \R_+$, statement (i) implies
$$
  \Pp_{\min}^{\alpha+\delta}=\Pp_{\min}^{\alpha}= 
 \Pp_{\max}^{\alpha} =\Pp_{\max}^{\alpha+\delta} = \emptyset.
$$
Thus, relations~\eqref{eq:PminPmaxEpsDeltaA} and \eqref{eq:PminPmaxEpsDeltaB} trivially hold for every $\alpha>\alpha^\ast$ and any choice of $\delta \in \R_+$.
We next show that, for any $\epsilon \in \R_+$, the choice
\begin{equation}
\label{eq:choiceOfDelta}
\delta = \epsilon\, \min \left\{ \alpha, \frac{\alpha^\ast- \alpha}{1+\epsilon} \right\}
\end{equation}
is such that~\eqref{eq:PminPmaxEpsDeltaA} and \eqref{eq:PminPmaxEpsDeltaB} hold if $\alpha<\alpha^\ast$. Obviously, \eqref{eq:choiceOfDelta} implies $ \delta <\alpha^\ast-\alpha$. Thus, both $\Rr_\infty^{\alpha}$ and $\Rr_\infty^{\alpha+\delta}$ are contained in $\X$ and we obtain $\Pp^{\alpha}_{\min}=\Rr_\infty^{\alpha}$ and 
$\Pp^{\alpha+\delta}_{\min}=\Rr_\infty^{\alpha+\delta}$ according to~\eqref{eq:Pmin}.
Moreover, \eqref{eq:choiceOfDelta} implies $\delta \leq \epsilon\,\alpha$. We hence find
\begin{equation}
\label{eq:RInftyAlphaDeltaSubset}
\Rr_\infty^{\alpha} \subseteq \Rr_\infty^{\alpha+\delta} = (\alpha+\delta) \Rr_\infty^{1} = \frac{\alpha+\delta}{\alpha} \Rr_\infty^{\alpha} \subseteq (1+\epsilon)  \Rr_\infty^{\alpha}
\end{equation}
according to Lem.~\ref{lem:RLimitAlpha}, which proves~\eqref{eq:PminPmaxEpsDeltaA}.
To see that~\eqref{eq:PminPmaxEpsDeltaB} holds as well, we set $\eta := (1+\epsilon)^{-1} \in (0,1)$ and note that $\delta \leq (1-\eta)(\alpha^\ast - \alpha)$. We thus find
$\eta \, \Pp_{\max}^\alpha \subseteq \Pp_{\max}^{\alpha+\delta} \subseteq \Pp_{\max}^\alpha$
according to Lem.~\ref{lem:SkAffectedByDelta},
 which confirms~\eqref{eq:PminPmaxEpsDeltaB}. So far, we showed that~\eqref{eq:PminPmaxEpsDeltaA} and~\eqref{eq:PminPmaxEpsDeltaB} hold if $\alpha \neq \alpha^\ast$. It remains to prove that~\eqref{eq:PminPmaxEpsDeltaA} and~\eqref{eq:PminPmaxEpsDeltaB} hold only if $\alpha \neq \alpha^\ast$. To this end, assume   $\alpha = \alpha^\ast$. Then  $\Pp_{\min}^\alpha$ and $\Pp_{\max}^\alpha$ are nonempty according to statement (ii). However, it is easy to see that $\Pp_{\min}^{\alpha+\delta}$ and $\Pp_{\max}^{\alpha+\delta}$ are empty for every choice of $\delta \in \R_+$. Thus, neither~\eqref{eq:PminPmaxEpsDeltaA} nor~\eqref{eq:PminPmaxEpsDeltaB} can hold for $\alpha =\alpha^\ast$. 

Statement (iv). This statement can be understood as a generalization of~\cite[Thm.~6.3]{Kolmanovsky1998}. As a consequence, the proof is inspired by the proof of that theorem.
We show that, whenever $\alpha \neq \alpha^\ast$, there always exists a $k \in \N$ such that 
\begin{equation}
\label{eq:SkPlus1EqualsSk}
\Ss_{k+1}^\alpha = \Ss_{k}^\alpha.
\end{equation}
From~\eqref{eq:sequenceSk} and Lem.~\ref{lem:PmaxPmin},
 \eqref{eq:SkPlus1EqualsSk} implies $\Pp_{\max}^\alpha = \Ss_{k}^\alpha$. We study~\eqref{eq:relationSkRkPlus1} to show that~\eqref{eq:SkPlus1EqualsSk} occurs after a finite number of iterations. 
We first address the case $\alpha < \alpha^\ast$ and define $\mu:=\left(1-\frac{\alpha}{\alpha^\ast} \right) \in (0,1)$. In this case,
we have 
\begin{equation}
\label{eq:XminusRkPlus1Superset}
 \X \ominus \Rr_{k+1}^\alpha \supseteq  \X \ominus \Rr_{\infty}^\alpha = \X \ominus \alpha \Rr_{\infty}^1 \supseteq  \mu\, \X
\end{equation}
for every $k \in \N$ by definition of $\Rr_\infty^1$, due to Lem.~\ref{lem:RLimitAlpha}, and according to \cite[Rem.~2.1]{Kolmanovsky1998} in combination with~\eqref{eq:criticalAlpha}. Since $A$ is strictly stable and since $\X$ is a C-set, there always exists a $L \in \N$ such that 
\begin{equation}
\label{eq:AjPlus1XSubset}
A^{L+1} \X \subseteq \mu \,\X.
\end{equation}
From \eqref{eq:XminusRkPlus1Superset}, \eqref{eq:AjPlus1XSubset}, and the fact that $\Ss_k^\alpha\subseteq \X$ for all $k\in\N$, we have 
$$
A^{L+1} \Ss_L^\alpha \subseteq  A^{L+1} \X \subseteq \mu\, \X \subseteq \X \ominus \Rr_{L+1}^\alpha.
$$
Clearly, this implies $\Ss_{L+1}^\alpha =  \Ss_L^\alpha$
 according to~\eqref{eq:relationSkRkPlus1}. 
 Thus, the choice $k=L$ is such that~\eqref{eq:SkPlus1EqualsSk} holds. It remains to prove~\eqref{eq:SkPlus1EqualsSk} for the case  
 $\alpha > \alpha^\ast$. To this end, we will show that there always exists an $L \in \N$ such that 
\begin{equation}
\label{eq:XMinusRL}
\X \ominus \Rr_{L}^\alpha = \emptyset.
\end{equation}
  Obviously, \eqref{eq:XMinusRL} implies $ \Ss_{L}^\alpha = \emptyset$ and consequently $ \Ss_{L+1}^\alpha = \emptyset$ according to~\eqref{eq:relationSkRkPlus1}.
  Thus, $k=L$ is such that~\eqref{eq:SkPlus1EqualsSk} holds.
  To find an $L$ that satisfies \eqref{eq:XMinusRL} first 
  note that, for every $\epsilon \in \R_+$, we can choose an $L \in \N$ such that $\Rr_\infty^{1} \subseteq (1+ \epsilon) \Rr_L^{1}$ (see Thm.~\ref{thm:outerApproxRInf1} below for details). 
We next choose such an $L$ for some $\epsilon<\frac{\alpha}{\alpha^\ast}-1$ and obtain
$$
\Rr_{L}^\alpha = \alpha \, \Rr_{L}^1 \supseteq \frac{\alpha}{1+\epsilon} \, \Rr_{\infty}^1 = \frac{\alpha}{\alpha^\ast\,(1+\epsilon)} \, \Rr_{\infty}^{\alpha^\ast}
$$ 
in accordance with Lem.~\ref{lem:RLimitAlpha}. Now, since $\frac{\alpha}{\alpha^\ast\,(1+\epsilon)}>1$, we find $\Rr_{L}^\alpha \nsubseteq \X$, which immediately implies~\eqref{eq:XMinusRL}. 

Statement (v). According to statement (iv),  there exists a finite $k \in \N$ such that  $\Pp_{\max}^\alpha= \Ss_{k}^\alpha$ if $\alpha<\alpha^\ast$. 
Moreover,  $\Pp_{\max}^\alpha= \Ss_{k}^\alpha$ is nonempty due to statement (ii) and thus a C-set according to Lem.~\ref{lem:SkAlphaSLimit}.

Statement (vi). 
By definition of $\alpha^\ast$ in~\eqref{eq:criticalAlpha}, 
it is easy to see that $\alpha  = \alpha^\ast$ implies $\partial \Rr_\infty^\alpha \cap \partial \X \neq \emptyset$
and thus $\partial \Pp_{\min}^\alpha  \cap \partial \X  \neq \emptyset $. To prove statement (vi).(a), it remains to show that 
\begin{equation}
\label{eq:partialPmiCapPartialX}
\partial \Pp_{\min}^\alpha  \cap \partial \X  = \emptyset 
\end{equation}
 whenever $\alpha \neq \alpha^\ast$.
Relation~\eqref{eq:partialPmiCapPartialX} obviously holds for every $\alpha>\alpha^\ast$ since $\Pp_{\min}^\alpha=\partial \Pp_{\min}^\alpha=\emptyset$ in this case. We prove that~\eqref{eq:partialPmiCapPartialX} also holds for every $\alpha \in (0,\alpha^\ast)$ by contradiction. To this end, consider any  $\alpha \in (0,\alpha^\ast)$, let $\mu:=\frac{\alpha^\ast}{\alpha}$, and assume~\eqref{eq:partialPmiCapPartialX} is violated, i.e.,  there exists an $x \in \partial \Pp_{\min}^\alpha  \cap \partial \X$. Since $\X$ is a C-set, we obviously have $0 \notin \partial \X$ and thus $x \neq 0$. According to Lem.~\ref{lem:RLimitAlpha}, we easily prove $x^\ast:=\mu \,x \in \partial \Rr_\infty^{\alpha^\ast}$. However, due to $\mu>1$, we also find $x^\ast \notin \X$ which contradicts $\partial \Rr_\infty^{\alpha^\ast} \subseteq \mathrm{cl}(\Rr_\infty^{\alpha^\ast}) \subseteq \mathrm{cl}(\X)=\X$. 

To prove statement (vi).(b), we first show that $\partial \Pp_{\min}^\alpha \cap \partial\Pp_{\max}^\alpha \neq \emptyset$ if $\alpha = \alpha^\ast$ by contradiction. To this end, note that $\alpha = \alpha^\ast$ implies $\Pp_{\min}^\alpha \neq \emptyset$ and $\Pp_{\max}^\alpha \neq \emptyset$ according to statement (i).
Since $\Ss_\infty^\alpha=\Pp_{\max}^\alpha$ is nonempty, $\Ss_\infty^\alpha$ is closed according to Lem.~\ref{lem:SkAlphaSLimit}. Taking Lem.~\ref{lem:PminPmaxEmpty} into account, we thus infer
\begin{equation}
\label{eq:closerPAlphaMin}
\partial \Pp_{\min}^\alpha \subseteq \mathrm{cl}(\Pp_{\min}^\alpha) \subseteq \mathrm{cl}(\Pp_{\max}^\alpha)= \Pp_{\max}^\alpha.
\end{equation}
Moreover, $\alpha = \alpha^\ast$ implies $\partial \Pp_{\min}^\alpha \cap \partial\X \neq \emptyset$ according to statement (vi).(a).
Now, consider any $x \in \partial \Pp_{\min}^\alpha  \cap \partial \X$, note that $x\neq 0$, and assume $\partial \Pp_{\min}^\alpha \cap \partial\Pp_{\max}^\alpha = \emptyset$. Clearly, since $x \in \partial \Pp_{\min}^\alpha$, the latter relation can only hold if
 $x \notin \partial\Pp_{\max}^\alpha$. According to~\eqref{eq:closerPAlphaMin}, this requires $x \in \mathrm{int}(\Pp_{\max}^\alpha)$. However, since $\X$ and $\Pp_{\max}^\alpha$ are both closed with $\Pp_{\max}^\alpha\subseteq \X$, we have $\mathrm{int}(\Pp_{\max}^\alpha) \cap \partial \X = \emptyset$ which contradicts $x \in \partial \X$.
It remains to show that
\begin{equation}
\label{eq:partialPminCapPartialPmax}
\partial \Pp_{\min}^\alpha  \cap \partial \Pp_{\max}^\alpha  = \emptyset 
\end{equation}
 whenever $\alpha \neq \alpha^\ast$. Not surprisingly, the proof is similar to the corresponding part in the proof of statement (vi).(a). Again, relation~\eqref{eq:partialPminCapPartialPmax} holds for every $\alpha>\alpha^\ast$ since $\Pp_{\min}^\alpha=\partial \Pp_{\min}^\alpha=\emptyset$ in this case. Moreover, to show that~\eqref{eq:partialPminCapPartialPmax} also holds for $\alpha<\alpha^\ast$, consider any  $\alpha \in (0,\alpha^\ast)$, let $\mu:=\frac{\alpha^\ast}{\alpha}$, and assume there exists an $x \in \partial \Pp_{\min}^\alpha  \cap \partial \Pp_{\max}^\alpha $ (i.e., \eqref{eq:partialPminCapPartialPmax} is violated). Since $\Pp_{\max}^\alpha$ is a C-set according to statement (v), we obviously have $0 \notin \partial \Pp_{\max}^\alpha$ and thus $x \neq 0$. 
Similar as above, we  find $x^\ast:=\mu \,x \in \partial \Rr_\infty^{\alpha^\ast}=\partial \Pp_{\min}^{\alpha^\ast}$ but $x^\ast \notin \Pp_{\max}^\alpha$. This  contradicts \eqref{eq:closerPAlphaMin} since $\Pp_{\max}^{\alpha^\ast} \subseteq \Pp_{\max}^{\alpha}$.

Statement (vii). In principle, the proof of statement (vi).(b) shows that $\emptyset \neq \partial \Pp_{\min}^\alpha \cap \partial\X \subseteq \partial \Pp_{\max}^\alpha$ in case that $\alpha=\alpha^\ast$. We thus have $\partial \Pp_{\max}^\alpha \cap \partial \X \neq \emptyset$ if $\alpha=\alpha^\ast$.
Now, consider any $x \in \Pp_{\max}^{\alpha^\ast} \cap \partial \X$ and note that $x \neq 0$. We show that $x \in \partial \Pp_{\max}^{\alpha}$ for every $\alpha <\alpha^\ast$ by contradiction. To this end, choose any $\alpha \in (0,\alpha^\ast)$, let $\mu:=\frac{\alpha^\ast}{\alpha}$, and  assume $x  \notin \partial \Pp_{\max}^{\alpha}$.  Since we have $\partial \Pp_{\max}^{\alpha^\ast} \subseteq \Pp_{\max}^{\alpha^\ast} \subseteq \Pp_{\max}^\alpha $, $x  \notin \partial \Pp_{\max}^{\alpha}$ implies $x \in \mathrm{int}(\Pp_{\max}^{\alpha})$. However, since $\X$ and $\Pp_{\max}^\alpha$ are both C-sets (by assumption and statement (v), respectively) with $\Pp_{\max}^\alpha\subseteq \X$, we find $\mathrm{int}(\Pp_{\max}^\alpha) \cap \partial \X = \emptyset$ which contradicts $x \in \partial \X$. It remains to show that $\partial \Pp_{\max}^\alpha \cap \partial \X = \emptyset$ whenever $\alpha>\alpha^\ast$. Clearly, this statement trivially holds since $\partial \Pp_{\max}^\alpha=\Pp_{\max}^\alpha=\emptyset$ if $\alpha>\alpha^\ast$. 
\end{proof}

Let us briefly discuss the consequences of  Thm.~\ref{thm:criticalScaling}. Clearly, statement (iii) expresses some continuity properties of the sets $\Pp_{\min}^\alpha$ and $\Pp_{\max}^\alpha$.
In this context, first note that the relations 
\begin{equation}
\label{eq:trivialInclusions}
\Pp_{\min}^{\alpha} \subseteq  \Pp_{\min}^{\alpha+\delta} \qquad \text{and} \qquad \Pp_{\max}^{\alpha+\delta} \subseteq  \Pp_{\max}^{\alpha}
\end{equation}
 trivially hold (for $\alpha < \alpha+\delta \leq \alpha^\ast$) since a larger disturbance set ($\D^\alpha \subset \D^{\alpha+\delta}$) increases reachability but decreases stabilizabilty. More interestingly, Eqs.~\eqref{eq:PminPmaxEpsDeltaA} and \eqref{eq:PminPmaxEpsDeltaB} state that the variations of the sets $\Pp_{\min}^{\alpha}$ and $\Pp_{\max}^{\alpha}$ resulting from a modified scaling $\alpha+\delta$ can be kept arbitrarily small by a suitable choice of $\delta$. To formalize this observation, we showed that, for every scaling factor $1+\epsilon>1$, there exists a $\delta>0$ such that the ``larger'' sets in~\eqref{eq:trivialInclusions} (i.e., $\Pp_{\min}^{\alpha+\delta}$ and $\Pp_{\max}^{\alpha}$) are contained in the scaled version of the corresponding ``smaller'' sets in~\eqref{eq:trivialInclusions}.
This continuity property holds, however, if and only if $\alpha \neq \alpha^\ast$. In other words, the sets $\Pp_{\min}^{\alpha}$ and $\Pp_{\max}^{\alpha}$ show a discontinuity around~$\alpha^\ast$. Clearly, this discontinuity is caused by the fact that 
$\Pp_{\min}^{\alpha}$ and $\Pp_{\max}^{\alpha}$ are nonempty for $\alpha \in (0,\alpha^\ast]$ but empty for $\alpha>\alpha^\ast$ (see statement (ii)).
As the following analysis shows, the transition to an empty set happens, however, abruptly and ``without warning''.
To see this, first note that $\Pp_{\min}^{\alpha}$ grows with $\alpha$ on the interval $(0,\alpha^\ast]$ (see Eq.~\eqref{eq:trivialInclusions} or Eq.~\eqref{eq:Pmin} in combination with  Lem.~\ref{lem:RLimitAlpha}). Consequently, the mRPI set takes its maximal size just before collapsing to the empty set (see also statement (vi)). 
 The variations of the MRPI set are slightly more intuitive. In fact, as apparent from~\eqref{eq:trivialInclusions}, the set $\Pp_{\max}^{\alpha}$ shrinks as $\alpha$ increases. However, since we have $\partial \Pp_{\max}^\alpha \cap \partial \X \neq \emptyset$ for every $\alpha \in (0,\alpha^\ast]$ according to statement (vii), the MRPI set still has contact to the boundary of the state constraints immediately before collapsing to the empty set. We will illustrate the transition of both $\Pp_{\min}^{\alpha}$ and $\Pp_{\max}^{\alpha}$ with some examples in Sect.~\ref{subsec:IllustrationProperties}.
Finally note that statements (iv) and (v) do not specify the ``structure'' of the MRPI set for the special case $\alpha = \alpha^\ast$. However, as discussed in Sect.~\ref{subsec:IllustrationProperties}, this ``gap'' is reasonable since $\Pp_{\max}^\alpha$ may or may not be finitely determined for $\alpha = \alpha^\ast$.

\subsection{A related study of controlled invariant sets}
\label{subsec:similarStudyOfCI}

As discussed above, statement (iii) in Thm.~\ref{thm:criticalScaling}
states that the sets $\Pp_{\min}^\alpha$ and $\Pp_{\max}^\alpha$ have a ``discontinuity'' around the CSF.
We show in the following that this property
is, in a sense, unique to RPI sets by proving that such a discontinuity
does not appear for the maximal controlled invariant (MCI) set.
 This observation is interesting since the MRPI set for \textit{autonomous} systems with  \textit{unknown} disturbances $d(k)$ offers many similarities to the MCI set for \textit{deterministic} systems with \textit{controllable} inputs. To clarify this similarity, assume for a moment that $d(k)$ in~\eqref{eq:system} is a control input and that $\D^\alpha$ in~\eqref{eq:constraints} describes input constraints.
Then, the following definition of the MCI set (which assumes $\X\subseteq \R^n$, $\D^\ast \subset \R^m$, and $\alpha \in \R_+$ as above) is indeed similar to Def.~\ref{def:MRPI}.

 \begin{defn}
 \label{def:MCI}
A set $\Cc \subseteq \R^n$ is called controlled invariant (CI) for system~\eqref{eq:system} with constraints~\eqref{eq:constraints} if (i)  $\Cc \subseteq \X$ and (ii) for every $x \in \Cc$, there exists a $d \in \D^\alpha$ such that $A\,x+E\,d \in \Cc$. The union of all CI sets for~\eqref{eq:system} with~\eqref{eq:constraints} is called the maximal controlled invariant (MCI) set for~\eqref{eq:system} with~\eqref{eq:constraints} and denoted by $\Cc_{\max}^\alpha$.
\end{defn}

Moreover, the MCI set can be characterized by the sequence
\begin{equation}
\label{eq:QSequence}
\Q_{k+1}^\alpha := A^{-1} (\Q_k^\alpha \oplus (- E \, \D^\alpha)) \cap \X \qquad \text{with} \qquad  \Q_{0}^\alpha:=\X,
\end{equation}
which is reminiscent of Eq.~\eqref{eq:sequenceSk}.
In fact, analogously to the findings in Lem.~\ref{lem:PmaxPmin}, the limit  $\Q_\infty^\alpha := \lim_{k \rightarrow \infty} \Q_k^\alpha = \bigcap_{k=0}^\infty \Q_k^\alpha$ equals the MCI set, i.e., $\Cc_{\max}^\alpha=\Q_\infty^\alpha$ (see, e.g., \cite[Thm.~3.1]{Blanchini1994}).
The sets $\Q_k^\alpha$ and the limit $\Q_\infty^\alpha$ are known to be C-sets in $\R^n$ if the following assumption holds in addition to Assum.~\ref{assum:matricesAndSets} (see, e.g., \cite[Props.~3.1 and 3.2, and Rem.~4.1]{Blanchini1994}). 

\begin{assum}
\label{assum:stability22}
The eigenvalues $\lambda \in \C$ of the matrix $A_{22}$ in Assum.~\ref{assum:matricesAndSets} are stable (i.e., $| \lambda|\leq 1$).
\end{assum}

Obviously, in contrast to Assum.~\ref{assum:stability}, the eigenvalues of the matrices $A_{11}$ and $A_{22}$ are not required to be strictly stable. We require, however, stability of the matrix $A_{22}$ since $\Q_\infty^\alpha$ will obviously be empty otherwise.
The following theorem establishes a statement on the continuity of the sets $\Q_\infty^\alpha$ as a function of $\alpha$.
In contrast to the corresponding statement (iii) in Thm.~\ref{thm:criticalScaling}, relation~\eqref{eq:MCIespDelta} holds for every $\alpha \in \R_+$ (and a suitable choice of $\delta \in \R_+$ depending on $\epsilon \in \R_+$). In other words, the MCI set does not show the discontinuity discovered for the MRPI (and the mRPI) set.

\begin{thm}
\label{thm:MCIespDelta}
Let Assums.~\ref{assum:matricesAndSets} and~\ref{assum:stability22} be satisfied. Then, for every $\epsilon\in \R_+$, there exists a $\delta \in \R_+$ such that
\begin{equation}
\label{eq:MCIespDelta}
\Cc_{\max}^\alpha \subseteq \Cc_{\max}^{\alpha+\delta} \subseteq (1+\epsilon) \,\Cc_{\max}^\alpha.
\end{equation}
\end{thm} 
 
 The proof of Thm.~\ref{thm:MCIespDelta} builds on the following lemma, which establishes the continuity of the sets $\Q_{k}^\alpha$.
 Note that Lem.~\ref{lem:QkAffectedByDelta} is similar to but less conservative than Lem.~\ref{lem:SkAffectedByDelta}.
 
\begin{lem}
\label{lem:QkAffectedByDelta}
Let Assums.~\ref{assum:matricesAndSets} and~\ref{assum:stability22} be satisfied and let $\delta \in \R_+$. Then
\begin{equation}
\label{eq:QkAffectedByDelta}
\frac{\alpha}{\alpha+\delta} \,\Q_{k}^{\alpha+\delta} \subseteq \Q_{k}^\alpha \subseteq \Q_{k}^{\alpha+\delta}
\end{equation}
 for every $k \in \N$.
\end{lem}

\begin{proof}
We prove the claim by induction. Relation~\eqref{eq:QkAffectedByDelta} obviously holds for
$k = 0$. It remains to show that~\eqref{eq:QkAffectedByDelta} implies
\begin{equation}
\label{eq:QkPlusOneAffectedByDelta}
\frac{\alpha}{\alpha+\delta} \,\Q_{k+1}^{\alpha+\delta} \subseteq \Q_{k+1}^\alpha \subseteq \Q_{k+1}^{\alpha+\delta}
\end{equation}
To this end, first note that we have
\begin{equation}
\label{eq:QkPlusOneBeta}
\Q_{k+1}^\beta 
 = A^{-1} (\Q_k^\beta \oplus (- \beta\,E \, \D^\ast)) \cap \X
 = \beta \left( A^{-1} \left( \frac{1}{\beta} \Q_k^\beta \oplus (- E \, \D^\ast) \right) \cap \frac{1}{\beta} \X \right)
\end{equation}
for all $\beta \in \R_+$ according to Eqs.~\eqref{eq:distributiveLaws} and~\eqref{eq:QSequence}. 
Moreover, we find
\begin{equation}
\label{eq:subsetQkPlusOne}
A^{-1} \left( \frac{1}{\alpha+\delta} \Q_k^{\alpha+\delta} \oplus (- E \, \D^\ast) \right) \cap \frac{1}{\alpha+\delta} \,\X \subseteq A^{-1} \left( \frac{1}{\alpha} \Q_k^\alpha \oplus (- E \, \D^\ast) \right) \cap \frac{1}{\alpha} \X
\end{equation}
since the relations
$$
\frac{1}{\alpha+\delta} \Q_k^{\alpha+\delta} \subseteq \frac{1}{\alpha} \Q_k^\alpha  \qquad \text{and} \qquad \frac{1}{\alpha+\delta} \X \subseteq \frac{1}{\alpha} \X
$$
hold due to~\eqref{eq:QkAffectedByDelta} and since $\X$ is a C-set by assumption. Now, evaluating~\eqref{eq:QkPlusOneBeta}
for the two cases $\beta=\alpha$ and $\beta = \alpha + \delta$, and taking~\eqref{eq:subsetQkPlusOne} into account leads to~\eqref{eq:QkPlusOneAffectedByDelta}.
\end{proof}

\begin{proof}[Proof of Thm.~\ref{thm:MCIespDelta}]
We clearly have $\Cc_{\max}^{\alpha} \subseteq \Cc_{\max}^{\alpha+\delta}$ due to $\D^\alpha \subset \D^{\alpha + \delta}$. Moreover, Lem.~\ref{lem:QkAffectedByDelta} implies
$$
 \Cc_{\max}^{\alpha+\delta} \subseteq \frac{\alpha+\delta}{\alpha}\,\Cc_{\max}^{\alpha} =  \left(1+\frac{\delta}{\alpha}\right) \Cc_{\max}^{\alpha}.
$$
Thus, the choice $\delta = \epsilon \alpha$ proves~\eqref{eq:MCIespDelta}.
\end{proof}

\section{Computation of the critical scaling factor}
\label{sec:algorithm}

All statements in the central Thm.~\ref{thm:criticalScaling} are related to the CSF $\alpha^\ast$ from~\eqref{eq:criticalAlpha}.
As a consequence, the application of the findings in Thm.~\ref{thm:criticalScaling} require the knowledge of (an approximation of) $\alpha^\ast$.
In this section, we thus address the numerical computation (or approximation) of the CSF. As a preparation, it is important to note that the optimization problem (OP) in~\eqref{eq:criticalAlpha} can be simplified based on Lem.~\ref{lem:RLimitAlpha}. In fact, the computation of $\alpha^\ast$ can be reduced to the identification of the largest scaling $\alpha$ such that $\alpha \,\Rr^1_\infty \subseteq \X $, i.e.,  $\alpha^\ast$ can be equivalently defined as 
\begin{equation}
\label{eq:criticalAlphaR1}
\alpha^\ast :=  \sup \,\{\alpha \in \R_+ \,|\, \alpha \,\Rr^1_\infty \subseteq \X \}.
\end{equation}
Nevertheless, the direct solution of the optimization problem (OP) in~\eqref{eq:criticalAlphaR1} is usually not tractable. First,  contrary to the MRPI set, the set $\Rr_\infty^1$ is generally not finitely determined and may be open. Second, even if there exists a finite $k$ such that $\Rr_k^1 = \Rr_\infty^1$, the explicit computation of $\Rr_k^1$ is demanding due to the involved Minkowski additions (see Eq.~\eqref{eq:sequenceRk}). 
We provide solutions to both problems in the following. To solve the first problem, we compute arbitrarily close inner and outer approximations of $\Rr_\infty^1$ that allow to derive arbitrarily precise over- and underestimations of $\alpha^\ast$. Regarding the second problem, it is essential to note that we are not interested in the computation (or approximation) of the set $\Rr_\infty^1$ itself but only the scaling $\alpha^\ast$. We can thus use  techniques introduced in \cite{Blanchini1999,Kolmanovsky1998,Rakovic2005} to solve an OP similar to~\eqref{eq:criticalAlphaR1} without explicitly  computing reachable sets. 

We start by deriving conditions that allow to compute over- and underestimations of $\alpha^\ast$. 
In principle, lower and upper bounds for $\alpha^\ast$ can be computed according to the trivial statements in Lems.~\ref{lem:alphaUpperBound} and \ref{lem:alphaLowerBound} (which we prove in the appendix).

\begin{lem}
\label{lem:alphaUpperBound}
Let Assums.~\ref{assum:matricesAndSets} and~\ref{assum:stability}  be satisfied, let $k \in \N$, and let  $\alpha^\ast$ be as in~\eqref{eq:criticalAlphaR1}. Define
\begin{equation}
\label{eq:alphaUb}
\ub{\alpha} :=  \max \,\{\alpha \in \R_+ \,|\, \alpha \, \Rr_k^1 \subseteq \X \},
%\ub{\alpha} :=  \max_{\mu} \,\mu \quad \text{s.t.} \quad \mu \, \Rr_k^1 \subseteq \X,
\end{equation}
Then, $\alpha^\ast \leq \ub{\alpha}$.
\end{lem}

\begin{lem}
\label{lem:alphaLowerBound}
Let Assums.~\ref{assum:matricesAndSets} and~\ref{assum:stability}  be satisfied, let $\epsilon \in \R_+$, and let  $\alpha^\ast$  be as in~\eqref{eq:criticalAlphaR1}. 
Assume $k\in \N$ is such that
\begin{equation}
\label{eq:epsApproximationRinfty}
\Rr_\infty^1 \subseteq (1+\epsilon)\,\Rr_k^1,
\end{equation}
and define $\ub{\alpha}$ as in~\eqref{eq:alphaUb}.
Then, $ \alpha^\ast  \geq (1+\epsilon)^{-1}\, \ub{\alpha}$.
\end{lem}

In some special cases, e.g., if $A$ is nilpotent, the overestimation $\ub{\alpha}$ from Lem.~\ref{lem:alphaUpperBound} can be used to exactly compute $\alpha^\ast$. This observation is summarized in  the following corollary, which builds on \cite[Rem.~4.2]{Kolmanovsky1998} and which is proven in the appendix.

\begin{cor}
\label{cor:exactAlphaAst}
Let Assums.~\ref{assum:matricesAndSets} and~\ref{assum:stability}  be satisfied, let $\eta \in [0,1)$, and let  $\alpha^\ast$  be as in~\eqref{eq:criticalAlphaR1}. 
Assume $k\in \N_+$ is such that $A^k = \eta \,I_n$
and define $\ub{\alpha}$ as in~\eqref{eq:alphaUb}.
Then, $\alpha^\ast  = (1-\eta)\,\ub{\alpha}$.
\end{cor}

In general, the condition in Cor.~\ref{cor:exactAlphaAst} will not be satisfied and we have to approximate $\alpha^\ast$ using Lems.~\ref{lem:alphaUpperBound} and \ref{lem:alphaLowerBound}.
Obviously, accurately approximating $\alpha^\ast$ requires to  solve~\eqref{eq:alphaUb} for a set $\Rr_k^1$ that satisfies~\eqref{eq:epsApproximationRinfty}.
As detailed in Sect.~\ref{subsec:numericalImplementation} below, the OP~\eqref{eq:alphaUb} can be efficiently solved using techniques from \cite{Blanchini1999,Kolmanovsky1998,Rakovic2005}.
It remains to identify a suitable $k \in \N$ such that~\eqref{eq:epsApproximationRinfty} holds for a given $\epsilon \in \R_+$. To this end, first note that similar problems were addressed in
\cite{Hirata2003,Rakovic2005}. 
Adapting the idea from \cite{Hirata2003}, we could (for increasing $k$) compute the maximal RPI for the system~\eqref{eq:system} with state and disturbance constraints $(1+\epsilon)\,\Rr_k^1$ and $\D^\ast$, respectively. In the case that the resulting maximal RPI is nonempty, we infer that \eqref{eq:epsApproximationRinfty} holds (see \cite{Hirata2003} for details). Clearly, the procedure is computationally demanding since we have to explicitly compute $\Rr_k^1$ and the corresponding maximal RPI for multiple values of~$k$.
A more efficient method can be found in \cite{Rakovic2005}. This method, however, requires the set $E \D^\ast$ to be full-dimensional in $\R^n$ (see \cite[beginning of Sect. II]{Rakovic2005}).
Clearly, this condition is not necessarily guaranteed by Assum.~\ref{assum:matricesAndSets}.
Nevertheless, as illustrated in the following, we can easily adapt the results from \cite{Rakovic2005}. First, it is easy to see that there always exists an $M \in \N_{[1,r]}$ such that
\begin{equation}
\label{eq:W}
\W:=\bigoplus_{k=0}^{M-1} A_{11}^k E_1\, \D^\ast
\end{equation}
is a C-set in $\R^r$. In fact, since $\D^\ast$ is a C-set in $\R^m$, $\W$ is C-set if (and only if) the matrix
\begin{equation}
\label{eq:controllabilityMatrix}
\Gamma=\begin{pmatrix}
A_{11}^0 \,E_1 & \dots & A_{11}^{M-1} E_1 
\end{pmatrix} \in \R^{r \times (M m)}
\end{equation} 
has full rank.  According to Assum.~\ref{assum:matricesAndSets}, $\Gamma$ is guaranteed to have  full rank for $M=r$ (but it may or may not have full rank for $M \in \N_{[1,r-1]}$).
This observation is used in the following lemma, which is inspired by \cite[Thm.~1]{Rakovic2005}.

\begin{lem}
Let Assums.~\ref{assum:matricesAndSets} and~\ref{assum:stability}  be satisfied, let $\eta \in (0,1)$, and let $M \in \N_{[1,r]}$ be such that $\W$ defined in~\eqref{eq:W} 
is a C-set in $\R^r$. 
Then, there exists an $N \in \N_{+}$ such
\begin{equation}
\label{eq:PhiJWSubsetEtaW}
A_{11}^{M N} \,\W\subseteq \eta \,\W.
\end{equation}
\end{lem}

The proof immediately follows from the facts that $\W$ is a C-set in $\R^r$ and that $A_{11}$ has strictly stable eigenvalues (by Assum.~\ref{assum:stability}).
According to the following theorem, the combination of (i) an $M$ such that $\W$ is a C-set in $\R^r$
and (ii) an $N$ satisfying~\eqref{eq:PhiJWSubsetEtaW} allows to compute a $k$ such that relation~\eqref{eq:epsApproximationRinfty} holds. 
To see this, note that~Eqs.~\eqref{eq:epsApproximationRinfty} and~\eqref{eq:outerApproxRInfty1} are equivalent for the choice $\eta = \epsilon \,(1+\epsilon)^{-1} \in (0,1)$.

\begin{thm}
\label{thm:outerApproxRInf1}
Let Assums.~\ref{assum:matricesAndSets} and~\ref{assum:stability}  be satisfied, let $\eta \in (0,1)$, and let $M \in \N_{[1,r]}$ be such that $\W$ defined in~\eqref{eq:W} 
is a C-set in $\R^r$. Assume $N \in \N$ is such that~\eqref{eq:PhiJWSubsetEtaW} holds and set $k=M N$.
Then, 
\begin{equation}
\label{eq:outerApproxRInfty1}
\Rr_\infty^1 \subseteq (1-\eta)^{-1}\, \Rr_{k}^1.
\end{equation}
\end{thm}

\begin{proof}
Let $\Lambda:=A_{11}^M$ and consider the sequence
\begin{equation}
\label{eq:sequenceRTilde}
\T_{j+1} := \Lambda \, \T_j \oplus \W \qquad \text{with} \qquad  \T_{0}:=\{0\}
\end{equation}
and its limit
$\T_\infty:=\lim_{j \rightarrow \infty} \T_j = \bigcup_{k=0}^\infty \T_j$.
Note that, in analogy to~Lem.~\ref{lem:RkAlphaRLimit}, $\T_\infty$ is bounded since $ \W$ is a C-set and since the eigenvalues of $\Lambda$ are strictly stable.
Moreover, according to \cite[Thm.~1]{Rakovic2005}, satisfaction of~\eqref{eq:PhiJWSubsetEtaW} implies that 
\begin{equation}
\label{eq:outerApproxT}
\T_\infty \subseteq (1-\eta)^{-1}\,\T_N. 
\end{equation}
Now, to prove~\eqref{eq:outerApproxRInfty1},
we identify a close link between the sequences defined by~\eqref{eq:sequenceRk} and~\eqref{eq:sequenceRTilde}.
%Now, to prove~\eqref{eq:outerApproxRInfty1}, we show that the sequences~\eqref{eq:sequenceRk} and~\eqref{eq:sequenceRTilde} are closely linked. 
As a preparation, note that the sets $\T_{j}$ are ``only'' $r$-dimensional
since~\eqref{eq:sequenceRTilde} only considers the reduced dynamics $(A_{11},E_1)$. We thus introduce the lifted sets
\begin{equation}
\label{eq:TPlusWPlus}
\T_j^+:=\left\{ \begin{pmatrix}
x_1 \\0
\end{pmatrix} \in \R^n \,\bigg|\, x_1 \in \T_j \right\}  \qquad \text{and} \qquad \W^+:=\left\{ \begin{pmatrix}
w \\0
\end{pmatrix} \in \R^n  \,\bigg|\, w \in \W \right\},
\end{equation}
which are $n$-dimensional by construction.
In the following, we show that 
\begin{equation}
\label{eq:relationRT}
\Rr^1_{M j}=\T^+_{j}
\end{equation}
 for every $j \in \N$.
To see this, first note that the relation holds for $j=0$ since $\Rr_0^1=\T^+_{0}=\{0\}$. It remains to show that~\eqref{eq:relationRT} implies $\Rr_{M (j+1)}^1=\T^+_{j+1}$.
This is easily proven as we obtain 
$$
\Rr_{M (j+1)}^1 \!=\! A^M \Rr_{M j}^1 \oplus \bigg(\bigoplus_{k=0}^{M-1} A^k E\, \D^\ast \! \bigg) \!=\!A^M \T_j^+ \!\oplus \W^+\! =\! \left\{ \begin{pmatrix}
\!x_1\!\! \\0
\end{pmatrix} \!\in \R^n \,\bigg|\, x_1 \!\in \Lambda \, \T_{j} \oplus \W^\ast \!\right\} \! =\!\T^+_{j+1}
$$
according to~\eqref{eq:sequenceRk}, due to~\eqref{eq:TPlusWPlus} and~\eqref{eq:relationRT}, by definition of $\Lambda$, and corresponding to Eqs.~\eqref{eq:sequenceRTilde} and~\eqref{eq:TPlusWPlus}. 
Clearly, \eqref{eq:relationRT} implies $\Rr_\infty^1=\T^+_\infty$. Thus, from~\eqref{eq:outerApproxT} in combination with~ \eqref{eq:relationRT}, we finally infer
$$
\Rr_\infty^1=\T^+_\infty \subseteq (1-\eta)^{-1}\,\T_N^+ = (1-\eta)^{-1}\,\Rr_{M N}^1 
$$
which proves~\eqref{eq:outerApproxRInfty1}. 
\end{proof}

Theorem~\ref{thm:outerApproxRInf1} in combination with Lems.~\ref{lem:alphaUpperBound} and \ref{lem:alphaLowerBound} suggests to use the following algorithm to accurately approximate $\alpha^\ast$.
In fact, as formalized in Thm.~\ref{thm:algorithm} further below, Alg.~\ref{alg:approximationOfCriticalScaling} allows the computation of lower and upper bounds on $\alpha^\ast$ that satisfy
\begin{equation}
\label{eq:accuracyBounds}
\lb{\alpha} \leq \alpha^\ast \leq \ub{\alpha} \qquad \text{and} \qquad \frac{\ub{\alpha}}{\lb{\alpha}}-1 = \epsilon
\end{equation}
for a given error bound $\epsilon \in \R_+$.

\begin{alg}
\label{alg:approximationOfCriticalScaling}
Approximation of the CSF $\alpha^\ast$ for any error bound $\epsilon \in \R_+$.
\vspace{-1mm}
\begin{enumerate}
\item[(i)] Choose the smallest $M \in \N_{[1,r]}$ such that $\Gamma$ in \eqref{eq:controllabilityMatrix} has full rank and define $\W$ as in~\eqref{eq:W}.
\item[(ii)] Set $\eta = \epsilon \,(1+\epsilon)^{-1}$ and choose the smallest $N \in \N_+$ such that~\eqref{eq:PhiJWSubsetEtaW} holds.
\item[(iii)] Set $k= M N$, compute $\ub{\alpha}$ according to~\eqref{eq:alphaUb}, and return bounds  $\lb{\alpha}=\ub{\alpha}\,(1+\epsilon)^{-1}$ and $\ub{\alpha}$. 
\end{enumerate}

\end{alg}

\begin{thm}
\label{thm:algorithm}
Let Assums.~\ref{assum:matricesAndSets} and~\ref{assum:stability}  be satisfied, let $\epsilon \in \R_+$, and let  $\alpha^\ast$  be as in~\eqref{eq:criticalAlphaR1}. Then, Alg.~\ref{alg:approximationOfCriticalScaling} computes  $\lb{\alpha}$ and $\ub{\alpha}$ such that \eqref{eq:accuracyBounds} holds.
\end{thm}

\begin{proof}
Step (i) of Alg.~\ref{alg:approximationOfCriticalScaling} guarantees that $\W$ in~\eqref{eq:W} is a C-set in $\R^r$. The choice of $\eta$ in step (ii) implies $(1-\eta)^{-1}=(1+\epsilon)$. Now, choosing $N$ such that~\eqref{eq:PhiJWSubsetEtaW} holds, yields
$$
\Rr_\infty^1 \subseteq (1+\epsilon)\, \Rr_{M N}^1
$$
according to Thm.~\ref{thm:outerApproxRInf1}.
Thus, computing $\ub{\alpha}$ as in~\eqref{eq:alphaUb} for $k=M N$
and setting $\lb{\alpha}=\ub{\alpha}\,(1+\epsilon)^{-1}$ implies
$\lb{\alpha} \leq \alpha^\ast \leq \ub{\alpha}$ corresponding to Lems.~\ref{lem:alphaUpperBound} and \ref{lem:alphaLowerBound}.
Equation~\eqref{eq:accuracyBounds} holds  by construction.
\end{proof}

The numerical implementation of Alg.~\ref{alg:approximationOfCriticalScaling} requires the solution of two non-trivial problems. In fact, depending on the shapes of the sets  $\X$ and $\W$, the choice of $N$ and the computation of $\ub{\alpha}$ in steps (ii) and (iii), respectively, may be computationally demanding. 
For polytopic sets  $\X$ and $\D^\ast$, however, Alg.~\ref{alg:approximationOfCriticalScaling} can be implemented efficiently  as described in the following section.

\subsection{Numerical implementation for polytopic constraints} 
\label{subsec:numericalImplementation}

Algorithm~\ref{alg:approximationOfCriticalScaling} requires the identification of the smallest $N$ such that~\eqref{eq:PhiJWSubsetEtaW} holds and the computation of $\ub{\alpha}$ as in~\eqref{eq:alphaUb} for $k=M N$ (where $M \in \N_{[1,r]}$ is such that $\Gamma$ in~\eqref{eq:controllabilityMatrix} has full rank). Both problems can efficiently be solved for polytopic sets $\X$ and $\D^\ast$, which can be written as
$$
\X=\{ x \in \R^n \,|\,H_x \,x \leq \boldsymbol{1}_{l_x} \} \qquad \text{and} \qquad  
\D^\ast=\{ d\in \R^m \,|\,H_d\, d \leq \boldsymbol{1}_{l_d} \}$$
for some $H_x \in \R^{l_x \times n}$, $H_d \in \R^{l_d \times m}$, and $l_x,l_d \in \N_+$. 
Clearly, if $\D^\ast$ is a polytope, the same holds for the set $\W$ in~\eqref{eq:W}, i.e., there exists a matrix $H_w \in \R^{l_w \times r}$ with $l_w \in \N_+$ such that $\W=\{ w\in \R^r \,|\,H_w\, w\leq \boldsymbol{1}_{l_w} \}$.
As a consequence, the support function
\begin{equation}
\label{eq:hW}
h_\W(v) := \sup_{w \in \W} v\, w 
\end{equation}
associated with $\W \subset \R^r$ and defined for row-vectors $v \in \R^{1\times r}$ can be evaluated by solving a linear program (LP). The support function provides the key to efficiently solve~\eqref{eq:PhiJWSubsetEtaW}. In fact, according to~\cite[Eq.~(10)]{Rakovic2005}, Eq.~\eqref{eq:PhiJWSubsetEtaW} holds if and only if 
\begin{equation}
\label{eq:hWForN}
\max_{i \in \N_{[1,l_w]}} h_{\W}(e_i^T H_w A_{11}^{M N}) \leq \eta, 
\end{equation}
where $e_i$ is the $i$-th unit vector in $\R^{l_w}$.
Obviously, for given $M,N \in \N_+$, condition \eqref{eq:hWForN} can be verified by solving $l_w$ LPs.
The efficient evaluation of~\eqref{eq:alphaUb} for $k=M N$ requires some preparation. 
In fact, instead of solving~\eqref{eq:alphaUb} directly, we solve a similar problem related to the sequence \eqref{eq:sequenceRTilde} introduced in the proof of Thm.~\ref{thm:outerApproxRInf1}. 
In this context, let
\begin{equation}
\label{eq:muUb}
\ub{\mu} := \max \,\{ \mu \in \R_+ \,|\, \mu \, \T_N \subseteq \X^{-} \},
\end{equation}
where $\X^{-}$ denotes an $r$-dimensional subset of $\X$ defined as
\begin{equation}
\label{eq:XMinus}
\X^{-}:= \left\{ x_1 \in \R^r \,\bigg|\, \begin{pmatrix}
x_1 \\ 0 
\end{pmatrix} \in \X \right\}.
\end{equation}
The following lemma shows that $\ub{\mu}$ and $\ub{\alpha}$ are indeed closely related.

\begin{lem}
\label{lem:alphaMuUb}
Let Assums.~\ref{assum:matricesAndSets} and~\ref{assum:stability}  be satisfied, let $M \in \N_{[1,r]}$ be such that $\W$ defined in~\eqref{eq:W} is a C-set in $\R^r$, assume $N \in \N$ is such that~\eqref{eq:PhiJWSubsetEtaW} holds, and set $k= M N$. Let $\ub{\alpha}$ and $\ub{\mu}$ be defined as in~\eqref{eq:alphaUb} and~\eqref{eq:muUb}, respectively. Then, $\ub{\alpha}=\ub{\mu}$.
\end{lem}

The proof of Lem.~\ref{lem:alphaMuUb} immediately follows from  relation~\eqref{eq:relationRT} in combination with~\eqref{eq:TPlusWPlus}.
The solution of~\eqref{eq:alphaUb} for $k=M N$ can thus be reduced to the solution of~\eqref{eq:muUb}. 
Now, it is easy to see that $\X^{-}$ from~\eqref{eq:XMinus} is also a polytope that can be written as $\X^{-} = \{  x_1 \in \R^r \,|\, H_x^{-} x_1 \leq \boldsymbol{1}_{l_x^{-}}\}$ based on some matrix $H_x^{-} \in \R^{l_x^{-} \times r}$ with $l_x^{-} \in \N_{[1,l_x]}$.
According to~\cite[Eq.~(12)]{Rakovic2005}, the OP in~\eqref{eq:muUb} can thus be solved as
\begin{equation}
\label{eq:computationOfMuUb}
\ub{\mu} = \min_{i \in \N_{[1,l_x^{-}]}} \frac{1}{\sum_{j=0}^{N-1} h_{\W}(\varepsilon_i^T H_x^{-} \Lambda^j) },
\end{equation}
where $\varepsilon_i$ is the $i$-th unit vector in $\R^{l_x^{-}}$ and where $\Lambda=A_{11}^M$ as in the proof of Thm.~\ref{thm:outerApproxRInf1}. Obviously, \eqref{eq:computationOfMuUb} can be evaluated by solving $N l_x^{-} \leq N l_x$ LPs.

\begin{rem}
 We showed that Alg.~\ref{alg:approximationOfCriticalScaling} can be efficiently implemented if $\X$ and $\D^\ast$ are polytopes. In fact, in this case, the computation of $N$ and $\ub{\alpha}=\ub{\mu}$ can be carried out based on~\eqref{eq:hWForN} and~\eqref{eq:computationOfMuUb}.
Both equations require the (multiple) evaluation of the support function $h_{\W}$ (see~\eqref{eq:hW}). As discussed above, $h_{\W}(v)$ can be evaluated by solving an LP.
Under certain conditions, $h_{\W}(v)$ can even be evaluated without solving an OP. In fact, if $\W$ can be described as an affine transformation of the hypercube in $\R^r$, $h_{\W}(v)$ can be computed analytically (see \cite[Rem.~3]{Rakovic2005}). 
However, it is in general not straightforward to link this condition to the underlying set $\D^\ast$. Nevertheless, an analytic evaluation of $h_{\W}(v)$ can even be guaranteed under less restrictive conditions. In fact, it is easy to see that $\W$ is a zonotope (i.e., the Minkowski sum of a finite number of line segments (see, e.g., \cite{Fukuda2004} for details)) whenever $\D^\ast$ is zonotopic.
Now, if $\W$ is a zonotope that can be written as
$
\W = \{ w \in \R^r \, |\, \exists \,\beta_1,\dots,\beta_L \in [-1,1]: w = \sum_{i=1}^L \beta_i \,z_i \}
$ 
for some $z_1,\dots,z_L \in \R^r$ with $L\in \N_+$, then the evaluation of $h_{\W}(v)$ results in $h_{\W}(v) = \sum_{i=1}^L | v \,z_i|$. This property is used to approximate $\alpha^\ast$ without any optimization for every example in Tab.~\ref{tab:resultsExamples} below.
\end{rem}

\section{Numerical examples}
\label{sec:examples}

In the following, we first illustrate the identified properties of RPI sets as summarized in Thm.~\ref{thm:criticalScaling} for three  simple illustrative examples. Afterwards, we discuss the related results on CI sets in Thm.~\ref{thm:MCIespDelta} with one example.
Finally, we apply Alg.~\ref{alg:approximationOfCriticalScaling} to approximate the CSF for a number of examples from the literature.

\subsection{Illustration of the identified properties of RPI and CI sets}
\label{subsec:IllustrationProperties}

\begin{exmp}
\label{exmp:Ex1}
Consider system~\eqref{eq:system} with $A=0.5$ and $E=1$ and the constraints $\X=[-2,2]$ and $\D^\ast=[-1,1]$. Note that the same system was also analyzed in \cite[Exmp.~6.10]{Blanchini2008} (without state constraints) and \cite[Exmp.~1]{SchulzeDarup2016_CDC}. 
For this simple example, the sets $\Ss_k^\alpha$ and $\Rr_k^\alpha$ of the sequences~\eqref{eq:sequenceSk} and~\eqref{eq:sequenceRk} can be stated explicitly. In fact, it is easy to prove that we have
\begin{align}
\label{eq:RkExample1}
\Rr_k^\alpha &=\{ x \in \R \,|\, |x| \leq \rho_k^\alpha \} &\text{with}& \qquad\rho_k^\alpha:=2\,\alpha \left(1 - 0.5^{k} \right) \qquad \text{and} \\
\label{eq:SkExample1}
\Ss_k^\alpha &= \{ x \in \R \,|\, |x| \leq \min \{2,\sigma_k^\alpha \} \}  &\text{with}& \qquad \sigma_k^\alpha:=2\,(\alpha -2^{k} (\alpha-1))
\end{align}
 for every $k \in \N$. 
Note that (in accordance with~\eqref{eq:relationSkRkPlus1}) $\rho_k^\alpha$ and $\sigma_k^\alpha$ are related to one another by  
$$
 \sigma_{k}^\alpha = (0.5^{k})^{-1} ( 2 - \rho_k^\alpha) = 2^{k}\, (2- \rho_k^\alpha).
$$ 
 Now, according to~\eqref{eq:RkExample1}, we have  $\lim_{k\rightarrow \infty} \rho_k^\alpha = 2\,\alpha$ but $\rho_k^\alpha  < 2\,\alpha$ for every $k \in \N$. We thus find $\Rr_\infty^\alpha=(-2\,\alpha,2\,\alpha)$ for every $\alpha \in \R_+$. Since $\Rr_\infty^\alpha \nsubseteq \X$ for every $\alpha>1$, we obtain $\alpha^\ast = 1$ according to~\eqref{eq:criticalAlpha}. Clearly, $\alpha^\ast$ is well-defined and finite as guaranteed by statement (i) in Thm.~\ref{thm:criticalScaling}. 
Regarding the limit $\Ss_{\infty}^\alpha$, for every $\alpha \leq 1$, we infer 
$\Ss_{\infty}^\alpha = \Ss_{1}^\alpha=\Ss_0^\alpha = [-2,2]$ from~\eqref{eq:SkExample1}. Moreover, for every $\alpha>1$, we have $\Ss_{k+1}^\alpha=\Ss_k^\alpha = \emptyset $ for every $k  \in \N$ with
\begin{equation}
\label{eq:boundOnKEx1}
k >\log_2 \left( \frac{\alpha}{\alpha-1} \right)
\end{equation}
and thus $\Ss_\infty^\alpha = \emptyset$. To see this, note that $\sigma_k^\alpha<0$ (and $\rho^\alpha_k>2$) for every $k \in \N$ satisfying~\eqref{eq:boundOnKEx1}. 
We thus obtain
\begin{equation}
\label{eq:PminPmaxEx1}
\Pp_{\min}^\alpha =  \left\{ \begin{array}{ll}
(-2\,\alpha,2\,\alpha) & \text{if} \quad \alpha \leq 1, \\
\emptyset & \text{otherwise},
\end{array}\right. \qquad \text{and} \qquad 
\Pp_{\max}^\alpha =  \left\{ \begin{array}{ll}
[-2,2] & \text{if} \quad \alpha \leq 1, \\
\emptyset & \text{otherwise}
\end{array}\right.
\end{equation}
in agreement with statement (ii) in Thm.~\ref{thm:algorithm}.
We next illustrate that the choice of $\delta$ as in~\eqref{eq:choiceOfDelta} is such that statement (iii) holds.
Consider, for example, $\alpha = 0.5$ and $\epsilon = 0.1$, then $\delta = 0.045 \leq 0.1 \min \{0.5, 0.\overline{45} \}$ is indeed such that 
\begin{align*}
\Pp_{\min}^{\alpha} &= (-1,1) \subseteq  \Pp_{\min}^{\alpha+\delta} = ( -1.09,1.09)\subseteq (1+\epsilon) \,\Pp_{\min}^{\alpha} = (-1.1,1.1) \quad \text{and} \\
\Pp_{\max}^{\alpha+\delta} &= [-2,2] \subseteq  \Pp_{\max}^{\alpha} = [ -2,2]\subseteq (1+\epsilon) \,\Pp_{\max}^{\alpha+\delta} = [-2.2,2.2].
\end{align*}
Regarding statement (iv), we showed above that $\Pp_{\max}^{\alpha}=\Ss_\infty^\alpha$ is finitely determined for every $\alpha \in \R_+$.
Equation \eqref{eq:PminPmaxEx1} also  implies that $\Pp_{\max}^\alpha $ is a C-set  in $\R^n=\R^1$ for every $\alpha \leq 1$ and thus confirms statement (v).
It is interesting to note that, for this example, $\Pp_{\max}^\alpha $ is finitely determined and a C-set even for $\alpha = \alpha^\ast$. Finally, statements (vi) and (vii) obviously hold since we have $\partial \Pp_{\min}^\alpha \cap \partial \X = \partial \Pp_{\min}^\alpha \cap \partial \Pp_{\max}^\alpha = \{-2,2\} \neq \emptyset$ if and only if $\alpha = \alpha^\ast = 1$, and $\partial \Pp_{\max}^\alpha \cap \partial \X = \{-2,2\} \neq \emptyset$ for every $\alpha \leq \alpha^\ast =1$.
\end{exmp}

Example~\ref{exmp:Ex1} confirms all statements in Thm.~\ref{thm:criticalScaling}. However, for this simple example, many statements (in particular statements (iii).(b), (v), and (vii)) are trivially fulfilled. We thus address another example to illustrate these findings for a slightly more complicated setup.

\begin{figure}[htp] 
\psfrag{a}[cb][c]{(a)}
\psfrag{b}[cb][c]{(b)}
\psfrag{c}[cb][c]{(c)}

\psfrag{x1}[ct][c]{$x_1$}
\psfrag{x2}[c][cb]{$x_2$}

\psfrag{x11}[ct][ct]{$-2$}
\psfrag{x12}[ct][ct]{$0$}
\psfrag{x13}[ct][ct]{$2$}

\psfrag{y11}[r][r]{$-2$}
\psfrag{y12}[r][r]{$0$}
\psfrag{y13}[r][r]{$2$}

\psfrag{x21}[ct][ct]{$-1$}
\psfrag{x22}[ct][ct]{$0$}
\psfrag{x23}[ct][ct]{$1$}

\psfrag{y21}[r][r]{$-1$}
\psfrag{y22}[r][r]{$0$}
\psfrag{y23}[r][r]{$1$}

\psfrag{x31}[ct][ct]{$-5$}
\psfrag{x32}[ct][ct]{$0$}
\psfrag{x33}[ct][ct]{$5$}

\psfrag{y31}[r][r]{$-5$}
\psfrag{y32}[r][r]{$0$}
\psfrag{y33}[r][r]{$5$}

\begin{flushright}
\includegraphics[trim=0mm 0mm 0mm 1mm, clip=true,width=1.0\linewidth]{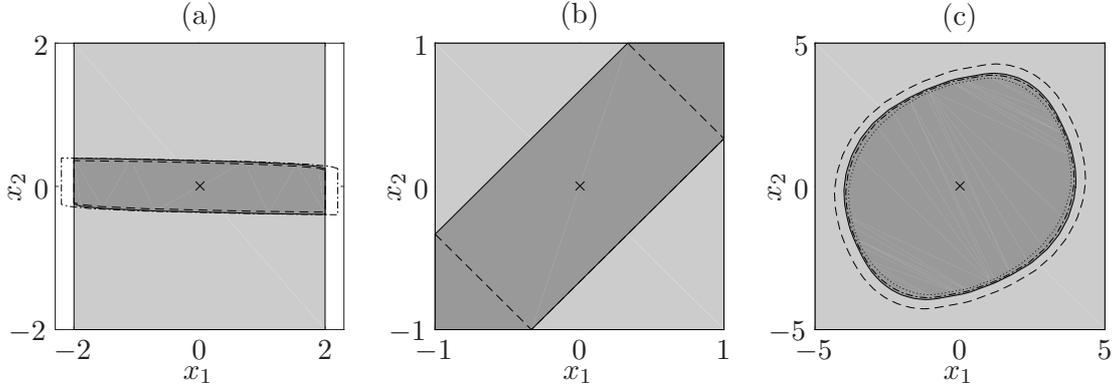}
\end{flushright}
\vspace{-4mm}
\caption{
Illustration of the sets in (a) Exmp.~\ref{exmp:Ex2}, (b) Exmp.~\ref{exmp:Ex3}, and (c) Exmp.~\ref{exmp:Ex4}.
In every figure, the light gray set visualizes the state constraints $\X$.
In addition, in (a), the sets $(1+\epsilon) \,\Pp_{\max}^{\alpha+\delta}$ (polytope with dash-dotted boundary),  $\Pp_{\max}^{\alpha}$ (dark gray polytope), and $\Pp_{\max}^{\alpha+\delta}$ (polytope with dashed boundary) are shown for $\alpha=0.5$, $\epsilon=0.1$, and $\delta = 0.045$. In (b), the sets $\Pp_{\min}^{\alpha}$ (polytope with dashed boundary) and $\Pp_{\max}^{\alpha}$ (dark gray polytope)  are illustrated for the special case $\alpha= \alpha^\ast = 0.\overline{6}$. In (c), the sets $\Q_{20}^{\alpha+\delta}$ (polytope with dashed boundary),  $\Q_{20}^{\alpha}$ (dark gray polytope), $\alpha\,(\alpha+\delta)^{-1} \,\Q_{20}^{\alpha+\delta}$ (polytope with dash-dotted boundary), and $\Q_{50}^{\alpha}$ (polytope with dotted boundary) are depicted for $\alpha=1$ and $\delta = 0.1$.} 
\label{fig:Exmps2to4}
\end{figure}

\begin{exmp}
\label{exmp:Ex2}
Consider system~\eqref{eq:system} with 
$$
A=\begin{pmatrix}
0.5 & 2.0 \\ 0.0 & 0.9
\end{pmatrix}
 \qquad \text{and} \qquad E=\begin{pmatrix}
1  \\ 0 
\end{pmatrix}
$$
 and the constraints $\X=\{ x\in \R^2 \,| \, \|x\|_\infty \leq 2 \}$ and $\D^\ast=[-1,1]$. 
Obviously, the matrices $A$ and $E$ offer the structure required in Assum.~\ref{assum:matricesAndSets} with $r=1$. 
 More interestingly, the matrices $A_{11}$ and $E_1$ are equivalent to the system matrices in Exmp.~\ref{exmp:Ex1}.
 Thus, the sets $\Rr_k^\alpha$ are given by
\begin{equation}
\label{eq:RkExample2}
 \Rr_k^\alpha= \{ x \in \R^2 \,|\, |x_1| \leq \rho_k^\alpha, \, x_2=0 \}
\end{equation}
with $\rho_k^\alpha$ as in~\eqref{eq:RkExample1} (see Rem.~\ref{rem:subsetDimr} or the proof of Thm.~\ref{thm:outerApproxRInf1} for details). 
An explicit description of the sets $\Ss_k^\alpha$ is more complicated but also not required for the following analysis.
However,  it is easy to see that the set $\X \ominus \Rr_k^\alpha$, which is involved in the computation of $\Ss_k^\alpha$ according to~\eqref{eq:relationSkRkPlus1}, can be written as
$$
\X \ominus \Rr_k^\alpha = \{ x \in \R^2 \,|\, |x_1| \leq 2-\rho_k^\alpha, \, |x_2| \leq 2 \}
$$
for every $k\in \N$. 
Now, analogously to Exmp.~\ref{exmp:Ex1}, we obtain $ \Rr_\infty^\alpha= \{ x \in \R^2 \,|\, |x_1| < 2\,\alpha, \, x_2=0 \}$ for every $\alpha \in \R_+$. Since we again have $\Rr_\infty^\alpha \nsubseteq \X$ for every $\alpha>1$, we again find $\alpha^\ast = 1$. While an explicit description of $\Ss_\infty^\alpha$ is not straightforward, it is easy to see that $\Ss_\infty^\alpha=\emptyset$ if (and only if) $\alpha>1$. In fact, for every $\alpha>1$, we have $\X \ominus \Rr_k^\alpha = \emptyset$ for every $k \in \N$ satisfying~\eqref{eq:boundOnKEx1} and thus $\Ss_\infty^\alpha=\Ss_{k+1}^\alpha=\Ss_k^\alpha = \emptyset$ according to~\eqref{eq:relationSkRkPlus1}. Statement (ii) in Thm.~\ref{thm:criticalScaling} thus also applies for this example.
Regarding statement (iii), consider again $\alpha=0.5$, $\epsilon = 0.1$, and $\delta = 0.045$ (satisfying condition~\eqref{eq:choiceOfDelta}). 
Clearly, statement (iii).(a) holds with the same reasoning as in Exmp.~\ref{exmp:Ex1}. In addition, the illustration of the sets 
$\Pp_{\max}^{\alpha+\delta}$,  $\Pp_{\max}^{\alpha}$, and $(1+\epsilon) \,\Pp_{\max}^{\alpha+\delta}$ in Fig.~\ref{fig:Exmps2to4}.(a) confirms statement (iii).(b).
We will only briefly address statements (iv) and (v). In fact, we only point out that, in contrast to Exmp.~\ref{exmp:Ex1}, the set $\Pp_{\max}^{\alpha}$ is not finitely determined and not a C-set for the special case $\alpha=\alpha^\ast=1$. 
To this end, we will show that $\Ss_\infty^1$ evaluates to
\begin{equation}
\label{eq:S1inftyEx2}
\Ss_\infty^1= \mathrm{cl}(\Rr_\infty^1)=  \{ x \in \R^2 \,|\, |x_1| \leq 2, \, x_2= 0 \}.
\end{equation}
Clearly, $\Ss_\infty^1$ as in~\eqref{eq:S1inftyEx2} is not a C-set in $\R^2$. Thus, since every set $\Ss_k^1$ is a C-set in $\R^2$ according to Lem.~\ref{lem:SkAlphaSLimit}, $\Ss_\infty^1$ and consequently $\Pp_{\max}^1$ cannot be finitely determined.
To see that~\eqref{eq:S1inftyEx2} holds, first note that we have 
$\mathrm{cl}(\Rr_\infty^1) = \mathrm{cl}(\Pp_{\min}^1) \subseteq \Pp_{\min}^1 = \Ss_\infty^1$ according to~\eqref{eq:closerPAlphaMin}. Now, $\mathrm{cl}(\Rr_\infty^1)\subset \Ss_\infty^1$ requires the existence of an $\xi \in \X$ with $\xi_2 \neq 0$ such that $\xi \in \Ss_\infty^1$. Such a $\xi$ does not, however, exist as we show next. 
As a preparation, consider the state  $\xi = (\,
0 \,\,\, \epsilon\,)^T$ for some $\epsilon \in \R_+$ and the disturbance sequence $d(j)=1 \in \D^1 = \D^\ast$ for every $j \in \N$. Since we have
$$
A^j = \begin{pmatrix}
0.5^j & 2 \, \sum_{i=0}^{j-1} 0.5^{j-1-i} \,0.9^i \\ 0.0 & 0.9^j
\end{pmatrix} = \begin{pmatrix}
0.5^j & 5 \, (0.9^j - 0.5^j) \\ 0.0 & 0.9^j
\end{pmatrix}
$$
for every $j\in \N$, we obtain
$$
x(k)= A^k \xi + \sum_{j=0}^{k-1} A^j E \,d(j) = \begin{pmatrix}
 5 \, (0.9^k - 0.5^k)\,\epsilon + \sum_{j=0}^{k-1} 0.5^j  \\  0.9^k\,\epsilon
\end{pmatrix} = \begin{pmatrix}
 2 +5\,\epsilon \, 0.9^k- (2+5\,\epsilon) \,0.5^k  \\  0.9^k\,\epsilon
\end{pmatrix}.
$$
Clearly, $x_1(k)>2$ and consequently $x(k) \notin \X$ for every $k \in \N$ with
$$
k > \log_{1.8} \left( \frac{0.4}{\epsilon} + 1 \right).
$$
In other words, for every $\epsilon \in \R_+$, there exists a finite $k \in \N$ such that $\xi \notin \Ss_k^1$ and thus $\xi \notin \Ss_\infty^1$. Analogously, one can show that $\xi = (\,
0 \,\,\, -\epsilon\,)^T\notin \Ss_\infty^1$ for any $\epsilon \in \R_+$ (by considering the disturbance sequence $d(j)=-1 \in  \D^\ast$ for every $j \in \N$). Thus, by convexity of $\Ss_\infty^1$, any $\xi \in \X$ with $\xi_2 \neq 0$ cannot be contained in $\Ss_\infty^1$, which proves~\eqref{eq:S1inftyEx2}. 
Finally, regarding statements (vi) and (vii) in Thm.~\ref{thm:criticalScaling}, first note that the state $\xi = (\,
2 \,\,\, 0\,)^T$ is contained in $\partial \X$ and $\partial \Ss_\infty^\alpha$ for every $\alpha \leq 1$. Moreover, $\xi \in \partial \Rr_\infty^1$. We thus have $\partial \Pp_{\max}^\alpha \cap \partial \X \neq \emptyset $ for every $\alpha \leq \alpha^\ast $ (statement (vii)) and   $\partial \Pp_{\min}^\alpha \cap \partial \X \neq \emptyset $ and $\partial \Pp_{\min}^\alpha \cap \partial \Pp_{\max}^\alpha \neq \emptyset $ if $\alpha = \alpha^\ast$ statement (vi)). It is easy to see that the two latter relations hold only if $\alpha=\alpha^\ast$, since $\Pp_{\max}^\alpha$ is a C-set in $\R^2$ with $\pm \xi \in \Pp_{\max}^\alpha$ for every $\alpha<1$.
\end{exmp}

Examples~\ref{exmp:Ex1} and \ref{exmp:Ex2} differ in that $\Pp_{\max}^{\alpha}$ for $\alpha =\alpha^\ast$ is finitely determined for the first but not for the second example.
However, Exmps.~\ref{exmp:Ex1} and \ref{exmp:Ex2} also offer some similarities. In fact, the set $\Rr_\infty^1$ is not closed for both examples, i.e., $\Rr_\infty^1 \neq \mathrm{cl}(\Rr_\infty^1)$.
Moreover, for $\alpha =\alpha^\ast$, the mRPI and MRPI sets are almost identical in the sense that $ \mathrm{cl}(\Pp_{\min}^\alpha) = \Pp_{\max}^\alpha$ for both examples. We analyze another example to point out that both relations do not hold in general and that there exists systems for which $\Rr_\infty^1 = \mathrm{cl}(\Rr_\infty^1)$
and $ \mathrm{cl}(\Pp_{\min}^{\alpha^\ast}) \neq \Pp_{\max}^{\alpha^\ast}$.

\begin{exmp}
\label{exmp:Ex3}
Consider system~\eqref{eq:system} with 
$$
A=\begin{pmatrix}
-0.5 & 0.5  \\ -0.5 & 0.5 
\end{pmatrix}
 \qquad \text{and} \qquad E=I_2
$$
 and the constraints $\X=\{ x\in \R^2 \,| \, \|x\|_\infty \leq 1 \}$ and $\D^\ast=\{ d\in \R^2 \,| \, \|d\|_1 \leq 1 \}$ 
as in \cite[p. 114]{Kouvaritakis2015}.
 Obviously, $A$ is nilpotent since $A^2=0$. 
 We consequently obtain $\Rr_1^1 = E \,\D^\ast = \D^\ast$ and
\begin{equation}
\label{eq:Rinfty1Ex3}
 \Rr_\infty^1 = \Rr_2^1 = A \,\D^\ast \oplus \D^\ast = 
\mathrm{conv} \left\{ \begin{pmatrix}
-1.5   \\ -0.5  
\end{pmatrix},\,\begin{pmatrix}
-0.5   \\ -1.5  
\end{pmatrix},\,\begin{pmatrix}
 1.5  \\0.5 
\end{pmatrix},\,\begin{pmatrix}
 0.5  \\1.5 
\end{pmatrix} \right\}
\end{equation}
according to~\eqref{eq:sequenceRk}. Obviously, $\Rr_\infty^1$ is closed, i.e., $\Rr_\infty^1 = \mathrm{cl}(\Rr_\infty^1)$. Now, from~\eqref{eq:Rinfty1Ex3} in combination with~\eqref{eq:criticalAlpha}, we infer $\alpha^\ast=0.\overline{6}$.
Regarding the computation of the sets $\Ss_k^\alpha$, we first find
$$
\Ss_1^\alpha = A^{-1} (\X \ominus \Rr_1^\alpha) \cap \X =  \{ x \in \R^2 \,|\, \|A \,x\|_\infty \leq 1-\alpha\} \cap \X= \{ x \in \X \,|\, |x_1-x_2| \leq 2 - 2\,\alpha \}
$$
according to~\eqref{eq:relationSkRkPlus1}. Moreover, due to $A^2=0$, we find
$$
(A^2)^{-1} (\X \ominus \Rr_2^\alpha ) = \left\{ \begin{array}{ll}
\R^2 & \text{if} \quad \alpha \leq \alpha^\ast \\
\emptyset & \text{otherwise}
\end{array} \right.
\quad \text{and thus} \quad
\Ss_\infty^\alpha = \Ss_2^\alpha = \left\{ \begin{array}{ll}
\Ss_1^\alpha & \text{if} \quad \alpha \leq \alpha^\ast \\
\emptyset & \text{otherwise}
\end{array} \right.
$$
For $\alpha=\alpha^\ast$, we consequently have
 $$\Pp_{\min}^\alpha  =
\mathrm{conv} \left\{ \begin{pmatrix}
-1.0   \\ -0.\overline{3}  
\end{pmatrix},\,\begin{pmatrix}
-0.\overline{3}   \\ -1.0  
\end{pmatrix},\,\begin{pmatrix}
 1.0  \\0.\overline{3} 
\end{pmatrix},\,\begin{pmatrix}
 0.\overline{3}  \\1.0 
\end{pmatrix} \right\}  \quad \text{and} \quad \Pp_{\max}^\alpha = \{ x \in \X \,|\, |x_1-x_2| \leq 0.\overline{6} \} $$
and thus $\Pp_{\min}^\alpha=\mathrm{cl}(\Pp_{\min}^\alpha) \subset \Pp_{\max}^\alpha$ as illustrated in Fig.~\ref{fig:Exmps2to4}.(b).
\end{exmp}

The main results of the paper are the identified properties of RPI sets as summarized in Thm.~\ref{thm:MCIespDelta}.
However, for ease of comparison, we also discussed some properties of CI sets in Sect.~\ref{subsec:similarStudyOfCI}.
We briefly illustrate the statement in Thm.~\ref{thm:MCIespDelta} or, more precisely, in Lem.~\ref{lem:QkAffectedByDelta} with the following example.

\begin{exmp}
\label{exmp:Ex4}
Consider system~\eqref{eq:system} with 
$$
A=\begin{pmatrix}
\blind{+}1.1 & 0.2 \\ -0.2 & 1.1
\end{pmatrix}
 \qquad \text{and} \qquad E=\begin{pmatrix}
0.5  & 0.0 \\ 0.0 & 0.2 
\end{pmatrix}
$$
 and the constrains $\X=\{ x\in \R^2 \,| \, \|x\|_\infty \leq 5 \}$ and $\D^\ast=\{ d\in \R^2 \,| \, \|d\|_\infty \leq 1 \}$ as in \cite[Exmp.~2.26]{SchulzeDarup2014_Phd}. 
Note that, analogously to  Sect.~\ref{subsec:similarStudyOfCI}, $d(k)$ in~\eqref{eq:system} and $\D^\alpha$ in~\eqref{eq:constraints} 
describe control inputs and input constraints for this example.
Further note that Assums.~\ref{assum:matricesAndSets} and \ref{assum:stability22} are satisfied (but not Assum.~\ref{assum:stability}). Computing the sets $\Q_k^\alpha$ and $\Q_k^{\alpha+\delta}$ according to~\eqref{eq:QSequence} for $\alpha=1$ and $\delta = 0.1$ yields the sets in Fig.~\ref{fig:Exmps2to4}.(c) for $k=20$. We clearly have $\alpha\,(\alpha+\delta)^{-1} \Q_{20}^{\alpha+\delta} \subseteq \Q_{20}^\alpha \subseteq \Q_{20}^{\alpha+\delta}$ as predicted by Lem.~\ref{lem:QkAffectedByDelta}. A graphical verification of Thm.~\ref{thm:MCIespDelta} is not possible for this example since $\Cc_{\max}^\alpha$ and $\Cc_{\max}^{\alpha+\delta}$ are not finitely determined.  However, the illustration of $\Q_{50}^\alpha$ in Fig.~\ref{fig:Exmps2to4}.(c) suggests that the sets $\Cc_{\max}^\alpha$ and $\Cc_{\max}^{\alpha+\delta}$ look similar to $\Q_{20}^\alpha$ and $\Q_{20}^{\alpha+\delta}$.
\end{exmp}

As discussed in Sect.~\ref{subsec:similarStudyOfCI}, Thm.~\ref{thm:MCIespDelta} shows that a critical scaling does not exist for CI sets associated with linear system and scaled input constraints. It thus points out an important difference between RPI and CI sets. Example~\ref{exmp:Ex4} illustrates another interesting difference. In fact, as apparent from Fig.~\ref{fig:Exmps2to4}.(c), the nonempty MCI set $\Cc_{\max}^\alpha \subseteq \Q_k^\alpha$ may have no contact points with the boundary of the state constraints while we have  $\partial \Pp_{\max}^\alpha \cap \partial \X \neq \emptyset$ for every nonempty MRPI set $\Pp_{\max}^\alpha$ according to statement (vii) in Thm.~\ref{thm:criticalScaling}.

\begin{landscape}
\begin{table}[htp]
\caption{\small Numerical results provided by Alg.~\ref{alg:approximationOfCriticalScaling} for the approximation of the CSF $\alpha^\ast$ for some examples.
For every example, the system matrices $A$ and $E$ and the constraints $\X$ and $\D^\ast$ are listed. Moreover, the dimension $r$, the numbers $M$ and $N$ computed in Alg.~\ref{alg:approximationOfCriticalScaling}, and the resulting bounds $\lb{\alpha}$ and $\ub{\alpha}$ are itemized.
For every example, the error bound $\epsilon$ was chosen as $\epsilon = 10^{-4}$. Note that, due to rounding errors, the second relation in~\eqref{eq:accuracyBounds} may not be exactly satisfied for the listed $\lb{\alpha}$ and $\ub{\alpha}$. For some examples, the state constraints depend on the (controller) matrix $K$. As detailed in the related references, we have    $K=(\,-0.6609 \,\,\, -1.3261 \, )$ for example 8, $K=(\,-0.89 \,\,\, -0.78 \, )$ for example 9, 
$K=(\,-0.77 \,\,\, -2.40 \,\,\, -2.59 \, )$ for example 10, and $K=(\,0.5484 \,\,\, 0.4299 \, )$  for example 11. For some examples, we slightly adjusted the state constraints so that $\X$ is a C-set as required in Assum.~\ref{assum:matricesAndSets}. In particular, the constraints $|x|\leq 2$ for example 1, $|x_1| \leq 2$ and $|x_2| \leq 1.5$ for example 4,   $\|x\|_\infty \leq 2$ for example 6, $|x_1| \leq 4$ and $|x_2| \leq 3$ for example 7,  $-50 \leq x_2$ for example 8, $-500 \leq x_1$ and $|x_2| \leq 800$ for example 10, and $|x_2|\leq 100$ for example 11 were not specified in the references.}
\label{tab:resultsExamples}
\centering
\footnotesize
\begin{tabular}{clccccccrcc}
\toprule 
no. & origin & $A$ & $E$ & $\X$ & $\D^\ast$ & $r$ & $M$ & $N$ & $\lb{\alpha}$ & $\ub{\alpha}$  \\
\midrule
1 & \cite[Exmp.~6.10]{Blanchini2008} & $0.5$ & $1$  & $|x|\leq 2$ & $|d| \leq 1$ & $1$ & $1$ & $14$ &  $0.999961$ & $1.000061$\\[3mm]
2 & here & $\begin{pmatrix}
0.5 & 2.0 \\ 0.0 & 0.9
\end{pmatrix}$ & $\begin{pmatrix}
1  \\ 0 
\end{pmatrix}$ & $\|x\|_\infty \leq 2 $ & $|d|\leq 1$ & $1$ & $1$ & $14$ &  $0.999961$ & $1.000061$\\[3mm]
3 & \cite[p. 114]{Kouvaritakis2015} &  
$\begin{pmatrix}
-0.5 & 0.5 \\ -0.5 & 0.5
\end{pmatrix}$  & $I_2$ & $\|x\|_\infty \leq 1$  & $\|d\|_1 \leq 1$  & $2$ & $1$ & $2$ & $0.666600$   & $0.666667$  \\[3mm]
4 & \cite[p. 200 f.]{Blanchini2008} &  
$\begin{pmatrix}
\blind{+}0.5 & 0.5 \\ -0.5 & 0.5
\end{pmatrix}$  & $\begin{pmatrix}
1  \\ 0 
\end{pmatrix}$ & $\begin{array}{rl}
|x_1|\!\!\!\! &\leq 2.0 \\
|x_2|\!\!\!\! &\leq 1.5
\end{array}$  & $|d| \leq 1$  &  $2$ & $2$ & $14$ & $0.857109$   & $0.857195$  \\[3mm]
%\cmidrule{1-10}
5 & \cite[Exmp.~4.1]{Kolmanovsky1995} &  $\begin{pmatrix}
0.5 & 0.0 \\ 1.0 & 0.1
\end{pmatrix}$  & $\begin{pmatrix}
1  \\ 1 
\end{pmatrix}$ & $\|x\|_\infty \leq 1$  & $ |d| \leq 1$ &  $2$ & $2$ & $8$ & $0.299977$   & $0.300007$ \\[3mm]
6 & \cite[Eq.~(13)]{Stoican2013} %\cite[Eq.~(17.15)]{Stoican2015} 
&  
$\begin{pmatrix}
0.9067 & -0.0687 \\ 0.0104 & \blind{+}0.7933
\end{pmatrix}$  & $\begin{pmatrix}
0.0272  \\ 0.3127 
\end{pmatrix}$ & $\|x\|_\infty \leq 2$  & $|d| \leq 1$ &  $2$  & $2$ & $60$ & $ 1.345374$   & $1.345509$  \\[3mm]
7 & \cite[Eq.~(17.16)]{Stoican2015} &  
$\begin{pmatrix}
-1.0559 & 1.1978 \\ -0.1711 & 0.9975
\end{pmatrix}$  & $\begin{pmatrix}
0.03  \\ 0.31 
\end{pmatrix}$ & $\begin{array}{rl}
|x_1|\!\!\!\! &\leq 4 \\
|x_2|\!\!\!\! &\leq 3
\end{array}$  & $|d| \leq 1$ &  $2$  & $2$ & $92$ & $0.992194$   & $0.992294$  \\[3mm]
8 & \cite[Sect.~4.1]{Mayne2005} &  
$\begin{pmatrix}
\blind{+}0.6696 & \blind{+}0.3369  \\ -0.6609 & -0.3261
\end{pmatrix}$  & $I_2$ & $\begin{array}{rl}
-50\leq x_2 \!\!\!\!&\leq 2 \\
|K x| \!\!\!\! & \leq 1
\end{array}$  & $\|d\|_\infty \leq 0.1$ &  $2$ & $1$ & $10$ & $3.362391$   & $3.362728$ \\[3mm]
9 & \cite[Exmp.~3.1]{Kouvaritakis2015} &  
$\begin{pmatrix}
-0.39 & -0.78 \\ \blind{+}0.50 & \blind{+}1.00
\end{pmatrix}$  & $-I_2$ & $\begin{array}{c}
-0.3 \leq x_1 \leq 0.7 \\
-0.5 \leq x_2 \leq 0.5 \\
-0.3 \leq K x \leq 0.2
\end{array}$  & $\|d\|_\infty \leq 0.05$ &  $2$ & $1$ & $21$ & $1.499907$   & $1.500057$ \\[4mm]
10 & \cite[Exmp.~4.2]{Kouvaritakis2015} &  
$\begin{pmatrix}
\blind{+}1.00 & \blind{+}1.00 & \blind{+}0.00 \\ \blind{+}0.00 & \blind{+}1.00 & \blind{+}1.00  \\ -0.77 & -2.40 & -1.59
\end{pmatrix}$  & $I_3$ & $\begin{array}{rl}
-500 \leq x_1\!\!\!\! &\leq 5 \\
  |x_2|\!\!\!\! &\leq 800 \\
 |K x|\!\!\!\! &\leq 4
\end{array}$  & $\|d\|_\infty \leq 0.25$ &  $3$ & $1$ & $15$ & $1.110404$   & $1.110515$ \\[5mm]
11 & \cite[Sect.~4]{Lee1999} &  
$\begin{pmatrix}
-0.2961 & -0.2300 \\ \blind{+}0.7058 & \blind{+}0.5500
\end{pmatrix}$  & $\begin{pmatrix}
1  \\ 1 
\end{pmatrix}$ & $\begin{array}{rl}
 |x_2|\!\!\!\! &\leq 100 \\
 |K x|\!\!\!\! &\leq 1
\end{array}$  & $ |d| \leq 0.7541$ &  $2$ & $2$ & $4$ & $1.007574$   & $1.007675$ \\
\bottomrule
\end{tabular}
\end{table}
\end{landscape}

\subsection{Approximation and interpretation of the critical scaling factor}
\label{subsec:ApproximationCSF}

We next approximate the CSF for some  examples from the literature. The underlying systems and the numerical results computed by Alg.~\ref{alg:approximationOfCriticalScaling} (using the procedures in Sect.~\ref{subsec:numericalImplementation}) are summarized in Tab.~\ref{tab:resultsExamples}. Note that the first three examples are identical to the first three examples in Sect.~\ref{subsec:IllustrationProperties}. Obviously, for these examples, the computed bounds $\lb{\alpha}$ and $\ub{\alpha}$ indeed satisfy $\lb{\alpha}\leq \alpha^\ast \leq \ub{\alpha}$ (with $\alpha^\ast$ as in Exmps.~\ref{exmp:Ex1}--\ref{exmp:Ex3}, respectively).
The computed bounds on $\alpha^\ast$ can also be easily verified for the fourth example in Tab.~\ref{tab:resultsExamples}.
In fact, for this example, we have $A^8=0.5^4 I_2 = 0.0625\,I_2$ so that Cor.~\ref{cor:exactAlphaAst} can be applied to exactly compute $\alpha^\ast$. 
Clearly, in analogy to Eq.~\eqref{eq:computationOfMuUb}, $\ub{\alpha}$ from Eq.~\eqref{eq:alphaUb} can be computed as
\begin{equation}
\label{eq:computationOfAlphaUb}
\ub{\alpha} = \min_{i \in \N_{[1,l_x]}} \frac{1}{\sum_{j=0}^{k-1} h_{\D}(\varepsilon_i^T H_x A^j E) }
\end{equation}
for polytopic sets $\X$ and $\D^\ast$ (where $\varepsilon_i$ now is the $i$-th unit vector in $\R^{l_x}$).
For $k=8$, we obtain $\ub{\alpha}=(1+0.5^3-0.5^5)^{-1}=1.09375^{-1}\approx 0.914286$ from~\eqref{eq:computationOfAlphaUb} and thus $\alpha^\ast = (1-0.0625)\,\ub{\alpha} = 0.9375 \,\ub{\alpha} \approx 0.857143$ according to Cor.~\ref{cor:exactAlphaAst}. Obviously, the bounds in Tab.~\ref{tab:resultsExamples} under- and overestimate $\alpha^\ast$ as expected. 

Having verified the results provided by Alg.~\ref{alg:approximationOfCriticalScaling} for four examples, we used the algorithm to compute CSFs for another seven  examples from the literature.
At this point, we have to comment on the usefulness and interpretation of the computed results. In fact, for most applications, the consideration of scaled disturbances as in~\eqref{eq:constraints} is not required.
An exception are parametric RPI sets as analyzed in \cite{Schaich2015,SchulzeDarup2016_CDC}, which can be used to describe state-dependent constraints (see \cite{Schaich2015} for details). However, even for conventional RPI sets, the computation of the CSF can be of interest. 
In fact, $\alpha^\ast$ can be understood as a measure for the actual robustness of RPI sets. To see this, note that constraints on the disturbances are usually not precisely known but rather estimations. 
Now, if the CSF for a system is smaller than $1$, the mRPI and MRPI sets for the system with the nominal disturbance constraints $\D^\ast$ are both empty and this would be recognized during the computation of $\Pp_{\min}^1$ or $\Pp_{\max}^1$. In contrast, if $\alpha^\ast\geq 1$, the computation of $\Pp_{\min}^1$ and $\Pp_{\max}^1$ results in nonempty sets and especially the nominal MRPI set contains no information about the ``closeness'' of the set to being empty. This information, however, can be easily inferred from $\alpha^\ast$. In fact, the closer $\alpha^\ast$ is to $1$ (from above) the closer $\Pp_{\min}^1$ and $\Pp_{\max}^1$ are to being empty and the less robust they are w.r.t.~uncertainties in the disturbance constraints $\D^\ast$. Based on this interpretation, we find that the nominal sets $\Pp_{\min}^1$ and $\Pp_{\max}^1$ for the two last examples in Tab.~\ref{tab:resultsExamples} are vulnerable to uncertainties of $\D^\ast$. Especially for the last example with $\alpha^\ast< 1.0077$, the sets $\Pp_{\min}^1$ and $\Pp_{\max}^1$ may be useless for any practical application. Note that this observation coincides with the analysis in  \cite[Sect.~4]{Lee1999}.
In fact, the authors of \cite{Lee1999} state that 
the maximal allowable bound for the disturbance of the fifth example is $|d|\leq 0.7841$. We used this constraint to define the nominal set $\D^\ast$ for this example. Obtaining a CSF close to $1$ thus confirms the results in \cite{Lee1999}. Moreover, the example shows that the idea of CSFs is foreshadowed in the literature (although it is not exactly specified in \cite{Lee1999} or elsewhere). Finally note that critical disturbances were also analyzed in \cite[Exmp.~6.3]{Kolmanovsky1998} for the fifth example in Tab.~\ref{tab:resultsExamples} (i.e.,  \cite[Exmp.~4.1]{Kolmanovsky1998}). The critical scaling $0.230769$ identified in \cite[Exmp.~6.3]{Kolmanovsky1998} does not coincide with our results in Tab.~\ref{tab:resultsExamples} since systems with additional output disturbances are addressed in \cite{Kolmanovsky1998} (see \cite[Eq.~(1.2)]{Kolmanovsky1998}). For the system class studied in this paper, it is, however, easy to see that the CSF evaluates to $\alpha^\ast=0.3$ for the fifth example (which coincides with the listed approximations).

\section{Conclusion}
\label{sec:conclusions}

This paper extends the theory on robust positively
invariant (RPI) sets for linear discrete-time systems
with additive disturbances. 
We presented a comprehensive analysis of the impact of scaled disturbance sets on the properties of the minimal and maximal RPI sets. In particular, we showed that there always exists a critical scaling factor (CSF), which determines the transition from nonempty to empty RPI sets. 
As summarized in Thm.~\ref{thm:criticalScaling}  - the main results of the paper - this CSF is  crucial for many properties of the mRPI and MRPI sets. Apart from the theoretical results in Thm.~\ref{thm:criticalScaling}, the computation of the CSF for a given system can be useful to quantify the robustness of RPI sets w.r.t.~uncertainties in the disturbance constraints (see Sect.~\ref{subsec:ApproximationCSF}). Moreover, knowledge of the CSF makes it possible to specify bounds on state and input constraints or acceptable magnitudes of disturbances when designing actuators, sensors,  and controllers.
To facilitate the application of the introduced analysis scheme, we 
derived an efficient algorithm for the approximation of the CSF $\alpha^\ast$ with arbitrary precision (see Alg.~\ref{alg:approximationOfCriticalScaling} and Thm.~\ref{thm:algorithm}). As summarized in Sect.~\ref{subsec:numericalImplementation}, the algorithm can be evaluated by solving a finite number of linear programs (LPs) if the constraints $\X$ and $\D^\ast$ are polytopes.

\section*{Acknowledgments}

Financial support by the German Research Foundation
(DFG) through the grants SCHU 2094/1-1 and SCHU 2094/2-1 is gratefully acknowledged.

\bibliographystyle{plain}

\begin{thebibliography}{10}


\bibitem{Blanchini1994}
F.~Blanchini.
\newblock {Ultimate boundedness control for uncertain discrete-time systems via
  set-induced Lyapunov functions}.
\newblock {\em IEEE Trans. Autom. Control}, 39(2):428--433, 1994.

\bibitem{Blanchini1999}
F.~Blanchini.
\newblock Set invariance in control.
\newblock {\em Automatica}, 35:1747--1767, 1999.

\bibitem{Blanchini2008}
F.~Blanchini and S.~Miani.
\newblock {\em Set-Theoretic Methods in Control}.
\newblock Birkh\"auser, 2008.

\bibitem{Fukuda2004}
K.~Fukuda.
\newblock From the zonotope construction to the {M}inkowski addition of convex
  polytopes.
\newblock {\em Journal of Symbolic Computation}, 38:1261--1272, 2004.

\bibitem{Glover1971}
J.~D. Glover and F.~C. Schweppe.
\newblock Control of linear dynamic systems with set constrained disturbances.
\newblock {\em IEEE Trans. Autom. Control}, 16(5):411--423, 1971.

\bibitem{Hirata2003}
K.~Hirata and Y.~Ohta.
\newblock {$\epsilon$}-feasible approximation of the state reachable set for
  discrete time systems.
\newblock In {\em Proc. of the 42nd IEEE Conference on Decision and Control},
  pages 5520--5525, 2003.

\bibitem{Kolmanovsky1995}
I.~Kolmanovsky and E.~G. Gilbert.
\newblock Maximal output admissible sets for discrete-time systems with
  disturbance inputs.
\newblock In {\em Proc. of American Control Conference}, pages 1995--1999,
  1995.

\bibitem{Kolmanovsky1998}
I.~Kolmanovsky and E.~G. Gilbert.
\newblock Theory and computation of disturbance invariance sets for
  discrete-time linear systems.
\newblock {\em Mathematical Problems in Engineering}, 4:317--367, 1998.

\bibitem{Kouvaritakis2015}
B.~Kouvaritakis and M.~Cannon.
\newblock {\em Model Predictive Control. {C}lassical, Robust and Stochastic.}
\newblock Springer, 2015.

\bibitem{Lee1999}
Y.~I. Lee and B.~Kouvaritakis.
\newblock Constrained receding horizon predictive control for systems with
  disturbances.
\newblock {\em International Journal of Control}, 72(11):1027--1032, 1999.

\bibitem{Mayne2006}
D.~Q. Mayne, S.~V. Rakovic, R.~Findeisen, and F.~Allgoewer.
\newblock {Robust output feedback model predictive control of constrained
  linear systems}.
\newblock {\em Automatica}, {42}({7}):{1217--1222}, {2006}.

\bibitem{Mayne2005}
D.~Q. Mayne, M.~M. Seron, and S.~V. Rakovi{\'c}.
\newblock Robust model predictive control of constrained linear systems with
  bounded disturbances.
\newblock {\em Automatica}, 41:219--224, 2005.

\bibitem{Rakovic2005}
S.~V. Rakovi{\'c}, E.~C. Kerrigan, K.~I. Kouramas, and D.~Q. Mayne.
\newblock Invariant approximations of the minimal robust positively invariant
  set.
\newblock {\em IEEE Trans. Autom. Control}, 50(3):406--410, 2005.

\bibitem{Rakovic2012}
S.~V. Rakovi{\'c}, B.~Kouvaritakis, M.~Cannon, and C.~Panos.
\newblock Fully parameterized tube model predictive control.
\newblock {\em International Journal of Robust and Nonlinear Control},
  22:1330--1361, 2012.

\bibitem{Schaich2015}
R.~M. Schaich and M.~Cannon.
\newblock Robust positively invariant sets for state dependent and scaled
  disturbances.
\newblock In {\em Proceedings of the 54th {C}onference on {D}ecision and
  {C}ontrol}, pages 7560--7565, 2015.

\bibitem{SchulzeDarup2014_Phd}
M.~Schulze~Darup.
\newblock {\em Numerical Methods for the Investigation of Stabilizability of
  Constrained Systems}.
\newblock PhD thesis, Ruhr-Universit{\"at} Bochum, 2014.

\bibitem{SchulzeDarup2016_CDC}
M.~Schulze~Darup, R.~M. Schaich, and M.~Cannon.
\newblock Parametric robust positively invariant sets for linear systems with
  scaled disturbances.
\newblock In {\em submitted to the 55th {C}onference on {D}ecision and
  {C}ontrol}, 2016.

\bibitem{Stoican2013}
F.~Stoican, M.~Hovd, and S.~Olaru.
\newblock Explicit invariant approximation of the {mRPI} set for {LTI} dynamics
  with zonotopic disturbances.
\newblock In {\em Proc. of 52th {C}onference on {D}ecision and {C}ontrol},
  pages 3237--3242, 2013.

\bibitem{Stoican2015}
F.~Stoican, C.~Oar{\u{a}}, and M.~Hovd.
\newblock {RPI} approximations of the {mRPI} set characterizing linear dynamics
  with zonotopic disturbances.
\newblock In S.~Olaru, A.~Grancharova, and F.~L. Pereira, editors, {\em
  Developments in Model-Based Optimization and Control}, pages 361--377.
  Springer, 2015.

\end{thebibliography}

\appendix

\section{Supplementary proofs}

\begin{proof}[Proof of Lem.~\ref{lem:SkRk}]
We prove the claim by induction. Relation~\eqref{eq:relationSkRk} obviously holds for $k = 0$ since we obtain $\Ss_0^\alpha = \X$ as in~\eqref{eq:sequenceSk}. It remains to show that~\eqref{eq:relationSkRk} implies
\begin{equation}
\label{eq:inductionStep}
\Ss_{k+1}^\alpha = \bigcap_{j=0}^{k+1} (A^{j})^{-1} (\X \ominus \Rr_{j}^\alpha).
\end{equation}
To this end, first note that~\eqref{eq:sequenceSk}  in combination with \eqref{eq:relationSkRk} yields
\begin{equation}
\label{eq:SkPlus1Rewritten}
\Ss_{k+1}^\alpha = A^{-1} \! \left( \left(\bigcap_{j=0}^k (A^{j})^{-1} (\X \ominus \Rr_{j}^\alpha) \right) \ominus E \,\D^\alpha\right) \cap \,\X\!
\end{equation}
The first term on the r.h.s.~in~\eqref{eq:SkPlus1Rewritten} can be rewritten as 
\begin{align}
\nonumber
A^{-1} \left( \left(\bigcap_{j=0}^k (A^{j})^{-1} (\X \ominus \Rr_{j}^\alpha) \right) \ominus E \,\D^\alpha\right) &= A^{-1} \left( \bigcap_{j=0}^k (A^{j})^{-1} \left((\X \ominus \Rr_{j}^\alpha) \ominus A^j  E \,\D^\alpha\right)  \right)  \\
\label{eq:termRewritten}
&=  \bigcap_{j=0}^k (A^{j+1})^{-1} \left(\X \ominus (\Rr_{j}^\alpha \oplus A^j  E \,\D^\alpha ) \right)  ,
\end{align}
where the rearrangements of the Pontryagin difference hold according to \cite[Thm.~2.1]{Kolmanovsky1998}.
Now, following the proof of \cite[Thm.~4.1]{Kolmanovsky1998}, 
the sequence in~\eqref{eq:sequenceRk} can be equivalently defined by $\Rr_{k+1}^\alpha = \Rr_k^\alpha \oplus A^k E \D^\alpha$. Using this relation in~\eqref{eq:termRewritten} and rewriting~\eqref{eq:SkPlus1Rewritten}, we finally obtain
\begin{equation}
\label{eq:finalLine}
\Ss_{k+1}^\alpha = \left( \bigcap_{j=0}^k (A^{j+1})^{-1} \left(\X \ominus \Rr_{j+1}^\alpha  \right)  \right) \cap \X = \left( \bigcap_{j=1}^{k+1} (A^{j})^{-1} \left(\X \ominus \Rr_{j}^\alpha  \right)  \right) \cap \X.
\end{equation}
Clearly, \eqref{eq:inductionStep} and~\eqref{eq:finalLine} are equivalent since $(A^{0})^{-1} (\X \ominus \Rr_{0}^\alpha)=I_n \, (\X \ominus \{0\})=\X$.
\end{proof}

\begin{proof}[Proof of Lem.~\ref{lem:RLimitAlpha}]
We first prove that $\Rr_k^\alpha = \alpha \,\Rr_k^1$ holds for every $k \in \N$ by induction. 
The relation obviously holds for $k=0$. Moreover, $\Rr_k^\alpha = \alpha \,\Rr_k^1$ implies $\Rr_{k+1}^\alpha = \alpha \,A \Rr_k^1 \oplus E \W^\alpha = \alpha \, \Rr_{k+1}^1$
according to Eqs.~\eqref{eq:distributiveLaws} and~\eqref{eq:sequenceRk}. We thus obtain $\Rr_\infty^\alpha = \alpha \,\Rr_\infty^1$ by definition of the limit $\Rr_\infty^\alpha$ in Lem.~\ref{lem:RkAlphaRLimit}. 
\end{proof}

\begin{proof}[Proof of Lem.~\ref{lem:alphaUpperBound}]
We obviously have $\Rr_k^1 \subseteq \Rr_\infty^1$ for every $k \in \N$. As a consequence, having $\alpha\, \Rr_\infty^1 \subseteq \X$ implies $\alpha \,\Rr_k^1 \subseteq \X$. We thus obtain  $\alpha^\ast\, \Rr_k^1 \subseteq \X$ and consequently $\alpha^\ast \leq \ub{\alpha}$.
\end{proof}

\begin{proof}[Proof of Lem.~\ref{lem:alphaLowerBound}]
Assume $\widehat{\alpha} \in \R_+$ is such that 
$\widehat{\alpha}\,(1+\epsilon)\,\Rr_k^1 \subseteq \X$. Then, \eqref{eq:epsApproximationRinfty} implies $\widehat{\alpha}\,\Rr_\infty^1 \subseteq \X$ and consequently $\widehat{\alpha} \leq \alpha^\ast$.
Now, from~\eqref{eq:alphaUb} it is easy to see that $\widehat{\alpha}= (1+\epsilon)^{-1} \,\ub{\alpha} $ is such that $\widehat{\alpha}\,(1+\epsilon)\,\Rr_k^1 \subseteq \X$, which proves the claim.
\end{proof}

\begin{proof}[Proof of Cor.~\ref{cor:exactAlphaAst}]
Following the argumentation in \cite[Rem.~4.2]{Kolmanovsky1998}, it is easy to see that $A^k = \eta \,I_n$ for some $k\in \N_+$ and $\eta \in [0,1)$ implies $\mathrm{cl}(\Rr_\infty^1) = (1-\mu)^{-1}\, \Rr_k^1$, which implies $(1-\mu)^{-1} \,\alpha^\ast  = \ub{\alpha}$.
\end{proof}

\end{document}